\documentclass[twocolumn]{autart}
\usepackage{graphicx}
\usepackage{mathtools}
\usepackage{txfonts}
\usepackage{fontawesome }
\usepackage{derivative}
\usepackage{booktabs}
\usepackage{empheq}
\usepackage{lscape}
 \usepackage{xcolor,colortbl}
\definecolor{lavender}{rgb}{0.9, 0.9, 0.98}
\usepackage{adjustbox}
\usepackage{caption}
\usepackage{subcaption}

\usepackage{tikz}
\usepackage{pgfplots}

\usepackage{enumitem}
\setitemize{nolistsep}
\setenumerate{nolistsep}
\usepackage{todonotes}
\usepackage{utfsym}

\newtheorem{example}{Example}

\allowdisplaybreaks

\DeclareMathOperator{\diag}{\ensuremath{\mathrm{diag}}}

\DeclareMathOperator{\Ortho}{O}

\newcommand{\ra}{\ensuremath{\rightarrow}}

\newcommand{\lra}{\ensuremath{\longrightarrow}}

\renewcommand{\ge}{\ensuremath{\geqslant}}

\renewcommand{\geq}{\ge}

\newcommand{\setmin}{\ensuremath{\setminus}}

\let\emptyset\varnothing

\newcommand{\abs}[1]{\ensuremath{\left\lvert{#1}\right\rvert}}

\newcommand{\secref}[1]{\S\ref{#1}}

\newcommand{\JOrtho}[1]{\Ortho_J}

\newcommand{\Borelsigalg}{\ensuremath{\mathfrak{B}}}

\newcommand{\inverse}{\ensuremath{^{-1}}}

\newcommand{\ol}{\overline}

\newcommand{\Let}{\coloneqq}
\newcommand{\nn}{\nonumber}

\newcommand{\Nz}{\ensuremath{\mathbb{N}}}

\newcommand{\norm}[1]{\ensuremath{\left\lVert #1 \right\rVert}}



\makeatletter
\renewcommand*\env@matrix[1][*\c@MaxMatrixCols c]{%
  \hskip -\arraycolsep
  \let\@ifnextchar\new@ifnextchar
  \array{#1}}
\makeatother

\DeclarePairedDelimiterX\aset[1]\lbrace\rbrace{\def\suchthat{\; \delimsize\vert\;}#1}
\DeclarePairedDelimiterX\LieBracket[2]\lbrack\rbrack{#1,#2}
\DeclarePairedDelimiterX\lcrc[2]\lbrack\rbrack{#1,#2}
\DeclarePairedDelimiterX\lcro[2]\lbrack\lbrack{#1,#2}
\DeclarePairedDelimiterX\lorc[2]\rbrack\rbrack{#1,#2}
\DeclarePairedDelimiterX\loro[2]\rbrack\lbrack{#1,#2}
\DeclarePairedDelimiterX\lm[2]\lparen\rparen{#1;#2}
\DeclarePairedDelimiterX\expecof[1]\lbrack\rbrack{#1}
\DeclarePairedDelimiterX\cexpecof[1]\lbrack\rbrack{\def\given{\: \delimsize\vert\:}#1}

\renewcommand{\theequation}{(\thesection.\arabic{equation})}
\makeatletter
\def\tagform@#1{\maketag@@@{\ignorespaces#1\unskip\@@italiccorr}}
\makeatother

\let\oldsqrt\sqrt
\def\sqrt{\mathpalette\DHLhksqrt}
\def\DHLhksqrt#1#2{%
\setbox0=\hbox{$#1\oldsqrt{#2\,}$}\dimen0=\ht0
\advance\dimen0-0.2\ht0
\setbox2=\hbox{\vrule height\ht0 depth -\dimen0}%
{\box0\lower0.4pt\box2}}

\newcommand{\Rbb}{\ensuremath{\mathbb{R}}}

\newcommand{\stt}[2]{\ensuremath{x^{#1}_{#2}}}
\newcommand{\cstt}[2]{\ensuremath{\breve{x}^{#1}_{#2}}}
\newcommand{\cont}[2]{\ensuremath{u^{#1}_{#2}}}
\newcommand{\ccont}[2]{\ensuremath{\breve{u}^{#1}_{#2}}}
\newcommand{\pmap}[1]{\ensuremath{g_{#1}}}
\newcommand{\cpmap}[1]{\ensuremath{\breve{g}_{#1}}}

\newcommand{\pe}[2]{\ensuremath{\zeta}^{#1}_{#2}}

\newcommand{\rnd}[2]{\ensuremath{f}_{#1}{#2}}

\newcommand{\pnoise}[2]{\ensuremath{w^{#1}_{#2}}}
\newcommand{\pnrv}[1]{\ensuremath{W_{#1}}}

\newcommand{\ps}{\ensuremath{\Omega_s}}

\newcommand{\PP}{\ensuremath{\mathsf{P}}}
\newcommand{\EE}{\ensuremath{\mathsf{E}}}

\newcommand{\Borel}{\ensuremath{\mathcal{B}}}

\newcommand{\node}{\ensuremath{N}}
\newcommand{\mnode}{\ensuremath{M}}
\newcommand{\grv}{\ensuremath{\mathcal{N}}}
\newcommand{\var}{\ensuremath{\mathsf{V}}}
\newcommand{\wvar}{\ensuremath{\mathsf{W}}}
\newcommand{\zv}[1]{\ensuremath{\mathbf{0}_{#1}}}
\newcommand{\pclass}[1]{\ensuremath{\Pi_{\mathsf{#1}}}}
\newcommand{\action}{\ensuremath{\ell}}

\newcommand{\caction}{\ensuremath{\breve{\ell}}}

\newcommand{\idist}{\ensuremath{\nu}}
\newcommand{\rate}{\ensuremath{r}}
\newcommand{\safe}{\ensuremath{Z}}
\newcommand{\unsafe}{\ensuremath{\breve{Z}}}
\newcommand{\meantr}[2]{\ensuremath{\mu}^{#1}_{#2}}
\newcommand{\eig}[2]{\ensuremath{\lambda}^{#1}_{#2}}
\newcommand{\ceig}[2]{\ensuremath{\breve{\lambda}}^{#1}_{#2}}
\newcommand{\ngb}{\ensuremath{\mathcal{N}}}
\newcommand{\dom}{\ensuremath{\mathbb{X}}}

\newcommand{\oprocendsymbol}{\hbox{$\bullet$}}
\newcommand{\oprocend}{\relax\ifmmode\else\unskip\hfill\fi\oprocendsymbol}

\newcommand{\stsp}{\mathbb{S}}
\newcommand{\sigalg}{\mathcal{G}}
\newcommand{\proc}{\ensuremath{X}}

\newcommand{\sample}{\omega_s}
\newcommand{\csample}{\breve{\omega}_s}
\newcommand{\policy}{\pi}
\newcommand{\cpolicy}{\breve{\pi}}
\newcommand{\control}{\ensuremath{U}}

\newcommand{\history}{\ensuremath{\ol{H}}}
\newcommand{\hist}{\ensuremath{\ol{h}}}
\newcommand{\phistory}{\ensuremath{H}}
\newcommand{\phist}{\ensuremath{h}}
\newcommand{\cphist}{\ensuremath{\breve{h}}}

\newcommand{\admact}{\mathbb{U}}
\newcommand{\QQ}{\mathsf{Q}}
\newcommand{\cQQ}{\breve{\mathsf{Q}}}

\newcommand{\sep}{\Sigma}
\newcommand{\abscont}{\ll}

\newcommand{\indic}[1]{\mathsf{1}_{#1}}

\newcommand{\tss}{\ensuremath{\Omega}}
\newcommand{\tsigalg}{\ensuremath{\mathcal{F}}}
\newcommand{\tsample}{\ensuremath{\omega}}
\newcommand{\tRR}{\ensuremath{\mathsf{R}}}

\newcommand{\tQQ}{\ensuremath{P}}
\newcommand{\ctQQ}{\ensuremath{\breve{P}}}

\DeclarePairedDelimiterX\probof[1]\lparen\rparen{#1}
\DeclarePairedDelimiterX\cprobof[1]\lparen\rparen{\def\given{\: \delimsize\vert\:}#1}

\begin{document}
\setlength{\abovedisplayskip}{2pt}
\setlength{\belowdisplayskip}{2pt}
\setlength{\abovedisplayshortskip}{2pt}
\setlength{\belowdisplayshortskip}{2pt}

\begin{frontmatter}

 \title{On the detection of the presence of malicious components in cyber-physical systems in the almost sure sense \thanksref{footnoteinfo}}

\thanks[footnoteinfo]{Souvik Das is supported by a PMRF grant from the Ministry of Human Resource Development, Govt. of India. Debasish Chatterjee acknowledges partial support of the SERB MATRICS grant MTR/2022/000656.}

\author[SD]{Souvik Das}\ead{souvikd@iitb.ac.in},
\author[PD]{Priyanka Dey}\ead{dey.p.aa@m.titech.ac.jp},
\author[SD]{Debasish Chatterjee}\ead{dchatter@iitb.ac.in}

\address[SD]{Systems \& Control Engineering, Indian Institute of Technology Bombay, Powai, Mumbai - 400076, India}
\address[PD]{Department of Systems and Control Engineering, Tokyo Institute of Technology, Tokyo 152-8550, Japan}



\begin{abstract}
	This article studies a fundamental problem of security of cyber-physical systems (CPSs): that of detecting, almost surely, the presence of malicious components in the CPS. We assume that some of the actuators may be \emph{malicious} while all sensors are \emph{honest}. We introduce a novel idea of \emph{separability} of state trajectories generated by CPSs in two situations: those under the nominal no-attack situation and those under the influence of an attacker. We establish its connection to security of CPSs in the context of detecting the presence of malicious actuators (if any) in them. As primary contributions we establish necessary and sufficient conditions for the aforementioned detection in CPSs modeled as Markov decision processes (MDPs). Moreover, we focus on the mechanism of perturbing the pre-determined control policies of the honest agents in CPSs modeled as stochastic linear systems, by injecting a certain class of random process called \emph{private excitation}; sufficient conditions for detectability and non-detectability of the presence of malicious actuators assuming that the policies are randomized history dependent and randomized Markovian, are established. Several technical aspects of our results are discussed extensively.
 \end{abstract}


\end{frontmatter}

\section{Introduction}
\label{sec:intro}
Cyber-Physical systems (CPSs) monitor and regulate several critical large-scale infrastructures such as smart grids, transportation systems, and wearable medical systems. Some recent cyber-attacks such as the Stuxnet computer worm attack \cite{ref:DK-13}, the cyber attack on Ukrainian's power grid \cite{ref:KZ-16}, and the Maroochy Shire water incident \cite{ref:NS-SM-17} demonstrate the security vulnerabilities of large-scale CPSs. With the increasing complexities of systems around us, the possibilities available to the attackers to launch \emph{sophisticated} and \emph{intelligent attacks} have increased. There is an emergent need to devote considerable attention to the issue of security of CPSs.  

This article studies the problem of detecting the presence of malicious attackers in CPSs.\footnote{A malicious attack is performed by an adversarial agent with the intent of degrading the performance of the underlying CPS, thereby restricting the CPS from achieving its goal.} In this article we focus on those attacks that affect the performance of the \emph{physical} layer of the CPSs consisting of sensors, actuators, and controllers connected over a network. 
Issues related to the \emph{cyber} layer of the CPS including cryptography, communication protocol, etc., that are associated with the underlying network of the CPS are not considered in this work.
 These attacks differ from the attacks on the \emph{cyber} layer of the system that typically involve several issues including cryptography, communication protocol, etc., associated with the underlying network of the CPS.

\textcolor{black}{Throughout this article, we stipulate that all sensors are honest while the malicious components may lurk in some of the actuators or the controllers and the communication channels between them.\footnote{A sensor is said to be \emph{honest} if it reports the true measurements that it observes. An actuator is \emph{honest} if it executes actions exactly in accordance with the prescribed control policy.} Motivated by the Stuxnet attack \cite{ref:DK-13} and the attack on Maroochy shire \cite{ref:NS-SM-17}, we consider \emph{arbitrarily intelligent attack} strategies employed by clever attackers. These attackers may collude with each other, share critical information and devise strategies to tamper with the nominal CPS, and remain undetected by acquiring \emph{complete} data about the states of the system at each instant of time by fully utilizing the underlying network of the system.
We emphasize that if all the actuators are malicious, then with the information on the system model and the states at each time, and with access to all the control inputs, the attackers can potentially hide forever; such attackers are not considered in this article. This marks a point of deviation from \emph{arbitrary} sensor attacks where the actuators and the controllers remain hidden from the adversaries and can consequently, take proper measures against them.}  

Several methods and important technologies have been developed over the last two decades to defend and secure a control system against malicious attacks; see \cite{ref:HS-VG-KHJ:22,ref:SH-MX-HHC-YL-14,ref:JG-DU-AC-etal-18,ref:DD-GLH-YX-etal-18,ref:HS-SA-KHJ-15,ref:AH-JL-FL-BL-17,ref:RM-IRC-14,ref:ST-DD-WZS-JY-SKD-16,ref:HH-JY:16} for surveys of such techniques and their applications. Among the most promising directions is the one that involves injecting \emph{private excitation} (also known as \emph{watermarks} in the literature) by the honest actuators into the system. They have been studied and used as active defense mechanisms against malicious attacks on the sensors and/or the actuators or both; see \cite{ref:KL-YM-KHJ-21} for a general overview. 

{\color{black}We start by reviewing the relevant works on watermarking techniques where the underlying system is considered to be linear. The term watermarking was first coined in the article \cite{ref:YM-BS-09} where it was used as an active defense mechanism against \emph{replay attacks} in the context of a linear time-invariant (LTI) system with Gaussian process noise. Under this premise, the technique attracted considerable attention over the past few years, and securing control systems by employing Gaussian-based watermarking against replay attacks (see \cite{ref:JRH-LDC-JGA-17,ref:CF-YQ-PC-WXZ-17,ref:LZ-KGV-JH-21} and \cite{ref:RG-CS-EN-SR-22} for co-design of watermarking signals and robust controllers), integrity attacks \cite{ref:YM-RC-BS-13,ref:SW-YM-BS-14,ref:YM-SW-BS-15}, stealthy attacks \cite{ref:YM-SW-BS-15,ref:MH-TT-VG-16}, and false data injection attacks in \cite{ref:YM-EG-AC-BS-10}, have been explored.

Security of LTI SISO/MIMO systems over a network under Gaussian random noise with complete and partial observations were considered in \cite{ref:BS-PRK-16,ref:BS-PRK-CDC-16,ref:BS-PRK-17}, wherein the authors assumed that the attack strategies employed by the adversary could be \emph{arbitrary}. A Gaussian-based dynamic watermarking scheme was proposed and the power distortion was used as a metric to detect the presence of malicious sensors. Based on the statistics of the output signals, several tests were introduced that ensure zero power distortion \emph{almost surely}. We highlight two notable features of these works: Firstly, the sensor attacks were considered to be arbitrary, and secondly, the reported result is concerned with the detection of sensor attacks \emph{almost surely}; these results are fundamentally stronger in the sense that it eliminates the possibility of error in detection. For SISO/MIMO LTI systems with partial observations, these results were generalized in \cite{ref:PH-MP-RV-AA-17} under a carefully designed attack model. The authors showed that persistent disturbances restrict the optimal design of the watermarking signal; to that end, they proposed an approach based on the internal model principle to compensate for persistent disturbance.  
An extension to the case of networked LTI systems was considered in \cite{ref:PH-MP-RV-AA-18} where the authors designed a watermarking signal based on null hypothesis testing. Under similar premises, several extensions were considered, e.g., to accommodate for the time-varying dynamics against replay attacks \cite{ref:MP-PH-AA-MJR-RV-20}, and for systems with non-Gaussian process noise \cite{ref:BS-PRK-19}. In addition, \cite{ref:RG-CS-EN-SR-22} addressed the co-design of watermarking signal and the authors in \cite{ref:SW-YM-BS-14,ref:HL-JY-YM-KHJ-18} introduced an optimization-based method for designing optimal watermarking signal.

More recently, a dynamic watermarking algorithm was proposed in \cite{ref:JT-JS-AG-21} for finite state-space and finite-action Markov decision processes where the authors obtained upper bounds on the mean time between false alarms, and the mean delay between the time an attack occurs and when it is detected. Various kinds of learning-based attacks were studied in \cite{ref:AF-WS-18,ref:MJK-AK-MF-TJ-19,ref:MJK-AK-MF-TJ-20}; see \cite{ref:JZ-LP-QLH-CC-SW-YX-21} for a survey on the applications of machine learning in attack detection in cyber-physical systems. We refer to Table \ref{tab:1} for a summary of all the related works in security of CPS where watermarking techniques are used.\footnote{We stress that by no means, this is an exhaustive list of contributions and apologies are extended for any omissions. We have only listed those articles that are relevant to our work.}

\begin{table*}[htbp]
\centering
\begin{tabular}{>{\columncolor{lavender}}ccccccc}
\toprule
\multicolumn{1}{c}{} & \multicolumn{2}{c}{\textbf{LTI systems}}&
\multicolumn{2}{c}{\textbf{LTV systems}}&
\multicolumn{1}{c}{\textbf{Nonlinear system}}&
\multicolumn{1}{c}{\textbf{MDP}}\\ 
\cmidrule(rl){2-3}\cmidrule(rl){4-5} 
\textbf{Attack types} & {G-WM} & {NG-WM} & {G-WM} & {NG-WM} & (NS)&
\\
\midrule
\textbf{Integrity} \\
Replay & \cite{ref:YM-BS-09,ref:YM-RC-BS-13} & \cite{ref:RG-CS-EN-SR-22,ref:JRH-LDC-JGA-17} & \cite{ref:MP-PH-AA-MJR-RV-20} &-&-&- \\
& \cite{ref:YM-SW-BS-15,ref:CF-YQ-PC-WXZ-17} & \cite{ref:LZ-KGV-JH-21,ref:SW-YM-BS-14} & \cite{ref:MP-SD-AJ-PH-AA-MJR-RV-20} &\\
& \cite{ref:HL-JY-YM-KHJ-18} & \cite{ref:YM-SW-BS-15,ref:PH-MP-RV-AA-20}& &\\
\textbf{Stealthy}  \\
Actuators & \cite{ref:MH-TT-VG-16} & - & -&-&-&-\\
Sensors & \cite{ref:YM-SW-BS-15} & -&- &-&-&- \\
\textbf{False data injection}  \\
Actuators &- &- &- &-&-&-\\
Sensors & \cite{ref:PH-MP-RV-AA-17,ref:MO-SS-PH-MP-RV-AA-20} & \cite{ref:PH-MP-RV-AA-20}& -&-&\cite{ref:MP-AJ-SD-QW-PH-AA-MJR-RV} & \cite{ref:JT-JS-AG-21} \\
\textbf{Arbitrary}  \\
Actuators  & \faCheck &- &\faCheck &-& \faCheck & \faCheckSquareO\\
Sensors &\cite{ref:BS-PRK-16,ref:BS-PRK-17} &\cite{ref:BS-PRK-16,ref:BS-PRK-19}  &- &- &\cite{ref:WHK-BS-PRK-19}&-\\
& \cite{ref:PH-MP-RV-AA-18,ref:TH-BS-PRK-LX-18} & & &\\
\textbf{Learning based}  \\
Actuators &\cite{ref:MJK-AK-MF-TJ-19} &\cite{ref:MJK-AK-MF-TJ-19} &- &-&\cite{ref:MJK-AK-MF-TJ-19} &-\\
Sensors &\cite{ref:MJK-AK-MF-TJ-19} &\cite{ref:MJK-AK-MF-TJ-19} &- &-&\cite{ref:MJK-AK-MF-TJ-19} &- \\
\bottomrule
\end{tabular}
\caption{This table summarizes the literature on the security of CPS where watermarking is adopted as a defense strategy and places our work in its appropriate niche. It provides a non-exhaustive list of works where the underlying systems are modeled as \emph{linear} or \emph{nonlinear systems}, and \emph{Markov decision processes} (MDP). The results of this article fit into this table at the location indicated by \faCheckSquareO  and \faCheck: \textcolor{black}{Theorem \ref{p:key} automatically applies to both linear and nonlinear systems (see \cite[\S1.2, p. 3]{ref:HLJL-12}, \cite[\S 3.5, p. 51]{ref:MLP-14}, \cite[Example 2.1]{ref:AraBorFerGhoMar-93} for more details).}
\\
\textsf{List of abbreviations used here}: LTI = Linear time-invariant; LTV = Linear time-varying; G-WM = Gaussian watermark; NG-WM = Non-Gaussian watermark.}
\label{tab:1}
\end{table*}

In most of the prior works, the sensors were considered to be vulnerable in the setting of linear CPSs model (except in \cite{ref:SW-YM-BS-14,ref:MH-TT-VG-16}). In contrast, our work generalizes the problem of arbitrarily intelligent \emph{actuator} attack detection to the setting of Markov decision processes (MDPs) and to the case of stochastic linear systems over a network. A common aspect of these preceding works is that the actuators were hidden from the adversaries, and dynamical watermarking techniques were particularly useful in those instances. On the other hand, we consider the actuators to be vulnerable, which contributes to the intricacies involved in designing a watermarking based control input for attack detection. }
We highlight the salient features of this work, below: 
\begin{itemize}[leftmargin=*,label = $\circ$]
    {\color{black} \item This article focuses on the class of attacks where a subset of the actuators is hijacked by adversaries.
    Types of well-known actuator attacks in the literature include denial-of-service, false data injection into the actuator channels, and eavesdropping. The actions of Stuxnet worms in the Iranian nuclear reactor \cite{ref:DK-13} is a classic example of actuator attacks and a sufficient indication that malicious attackers can and will \emph{arbitrarily} tamper with actuators/controller without being detected.

    \item Our main results --- Theorem \ref{p:key}, Theorem \ref{prop:hmain}, and Corollary \ref{prop:main}, stated in \S\ref{sec:separability} and \S\ref{sec:mrg}, \emph{remain agnostic to the specific functional nature of the control policies. In particular, we do not enforce any narrow limitations on the class of control policies adopted by the malicious actuators under the attack situation}; see \S \ref{subsec:types} for an elaborate discussion on the \emph{arbitrarily clever nature} of the malicious attackers, assumed herein.

    \item The results established in this article are related to \emph{almost sure} detection of the presence of malicious actuators. Such results are technically different from those that cater to assertions of the type `on average' or `with high probability', neither of which preclude the possibility of errors. }
\end{itemize}

\subsection*{Our contributions}
This article establish a general framework for \emph{arbitrarily intelligent actuator attack} detection. The key contributions of this article are as follows:
\begin{enumerate}[leftmargin=*]
    \item We introduce an idea of \emph{separability} of state trajectories between two classes of state trajectories generated by the system when it is influenced by all honest actuators versus those generated by the same system influenced by at least one malicious actuator.\footnote{A device that executes the task of separation is a \emph{separator}.} \textcolor{black}{We define the relationship between separability and security of CPSs in the context of detecting the presence of malicious activities in the system as evidenced via its trajectories (see Definition \ref{def:separator} and Proposition \ref{prop:det_from_sep} for further technical details) and establish the first abstract result (Theorem \ref{p:key}) in the context of CPSs modeled as Markov decision processes (MDPs).}  \textcolor{black}{It applies naturally to both linear and nonlinear stochastic systems by virtue of \cite[\S1.2, p. 3]{ref:HLJL-12}, \cite[\S 3.5, p. 51]{ref:MLP-14}, \cite[Example 2.1]{ref:AraBorFerGhoMar-93}.}

    \item For a linear time-invariant stochastic CPS with Gaussian process noise, we use the mechanism of injection of a random process called \emph{private excitation} with suitable statistics, by the honest actuators into the linear system for separating the state trajectories. In this article we choose the private excitation to be Gaussian.
    \begin{enumerate}
        \item Under the above premise, we provide sufficient conditions for the existence/non-existence of a \emph{separator} when the honest and corrupt policies are chosen as randomized and history dependent; see Theorem \ref{prop:main} for more details on these conditions.\footnote{We refer the readers to Definition \ref{def:policy} in Appendix \ref{appen:prob basics} for the definition of a randomized policy.}
        
        \item We also establish sufficient conditions for the existence and non-existence of a separator under the assumption that the honest and the corrupt policies are randomized and Markovian; see Corollary \ref{prop:hmain} in \S\ref{subsubsec:rh control} for more information on these conditions.
    \end{enumerate}
    The results are extensively discussed across several remarks in \S\ref{subsec:tech_discuss} that provide insights into the technical aspects of the notion of \emph{separability} and the established results.
 
    One of the key attributes of our results is that they are existential and rely on the statistics of infinitely long state trajectories; these conditions may not be verifiable in practice as they stand now. Moreover, our main results --- Theorem \ref{p:key}, Theorem \ref{prop:main} (and Corollary \ref{prop:hmain}) --- neither provide a mechanism for the detection of the specific malicious components nor do they provide an estimate of the number of such malicious components present therein. Notwithstanding their per se inapplicability, we regard our main results as \emph{baseline fundamental steps towards understanding the nature of the problem under consideration and assessing the boundary of possibilities}. The topic of implementable algorithm design for the detection of malicious components will be addressed in subsequent articles.
\end{enumerate} 

\subsection*{Organization}    
\textcolor{black}{\secref{sec:separability} formally introduces the idea of separability --- Definition \ref{def:separator} in \S\ref{subsec:concept_sep}, and discusses its connection to detecting the presence of malicious actuators in the CPS. The primary result of this article --- Theorem 1, is established in \S\ref{subsec:chief_obs} in the context of CPSs modeled as MDPs. In \S\ref{sec:mrg} we focus our attention on linear CPSs and give sufficient conditions under which the existence and non-existence of a separator can be asserted --- Theorem \ref{prop:main} and Corollary \ref{prop:hmain}. The types of attacks considered in this work are detailed in \S\ref{subsec:types}. Several technical examples and remarks in \S\ref{subsec:examples} and \S\ref{subsec:tech_discuss}, respectively, are provided to discuss the important ideas related to the main results in \S\ref{sec:separability} (Theorem \ref{p:key}) and  \S\ref{sec:mrg} (Theorem \ref{prop:main} and Corollary \ref{prop:hmain}). Finally, we conclude by giving a brief summary of our results and possible future directions.}

\subsection*{Notations}
We employ standard notation: The set of real numbers and positive integers are denoted by  \(\Rbb\) and \(\mathbb{N}\) respectively.
Let \(E\) be a topological space, then \(\Borelsigalg{(E)}\) denotes the Borel \(\sigma \mbox{-}\)generated by the topology of \(E\). \textcolor{black}{Consider two probability measures \(\PP, \PP'\) on the same measurable space \((\Omega, \sigalg)\). If \(\PP(A) = 0\) for every \(A\in\sigalg\) satisfying \(\PP'(A) = 0\), then \(\PP\) is \emph{absolutely continuous} \cite[p.\ 233, p.\ 438]{ref:Shi-16} with respect to \(\PP'\) and is denoted by \(\PP\abscont\PP'\).} For any square matrix \(A \), its determinant is denoted by \(\det(A)\). \textcolor{black}{By \(A \succ 0\) we mean that \(A\) is a positive definite matrix.} For any set \(S,\) the set \(S^{\mathsf{c}}\) denotes its complement. Given a set \(\dom\), we define the indicator function on a set \(Y \subset \dom\) by
\[
\dom \ni x \mapsto \indic{Y}(x) \Let \begin{cases}
1 \qquad \text{if }x \in Y,\\
0 \qquad \text{else.}
\end{cases}
\]

\section{Separability of state trajectories in a general context}
\label{sec:separability}
We introduce an idea of \emph{separability} of sample paths generated by a cyber-physical system (CPS) modeled as MDPs \cite{ref:PV-13,ref:JT-JS-AG-21} and establish its connection to the security of the CPS. 
\subsection{The setting of Markov decision processes (MDPs)}\label{subsec:MDP_CPS_Def}
{\color{black}
An MDP \cite{ref:AraBorFerGhoMar-93,ref:HLJL-12} with \(\stsp\) and \(\admact\) being the state-space and the admissible action space, respectively, is considered. Throughout this article, we stipulate that \(\stsp\) and \(\admact\) are nonempty Borel measurable subsets of the Euclidean space. 

We briefly review the framework of MDPs below: Define the \emph{canonical sample space} by \(\tss \Let (\stsp \times \admact)^{\Nz}\). Let \(\Borelsigalg{\bigl(\tss\bigr)}\) denotes the Borel \(\sigma\mbox{-}\)algebra induced by the product topology on \(\tss\) in a standard fashion.\footnote{See \cite[Sect. 2.1]{ref:AraBorFerGhoMar-93} for more information on the construction of the filtered measurable space.} An element \(\tsample \in \tss\) is a sequence of the form \(\tsample \Let \bigl(\stt{}{0},\cont{}{0},\stt{}{1},\cont{}{1},\ldots\bigr)\) where \(\stt{}{t}\in\stsp\) and \(u_t\in\admact\) are the state and the control action picked according to some policy \(\policy\) (see Definition \ref{def:policy} in Appendix \ref{appen:prob basics}) at time \(t\in\Nz\).\footnote{Here the class of all policies is denoted by \(\pclass{}.\) We refer the readers to \cite{ref:HLJL-12} for more details on different sub-classes of \(\pclass{}\).} The key processes of relevance are:
\begin{itemize}[leftmargin=*]
    \item the \emph{state process} \(\proc \Let (\proc_t)_{t \in \Nz}\) is defined by \(\proc_t (\tsample) \Let \stt{}{t}\);

    \item the \emph{control process} \(\control \Let (\control_t)_{t \in \Nz}\) is defined by \(\control_t (\tsample) \Let \cont{}{t}\);

    \item the \emph{history process} \(\history = (\history_t)_{t \in \Nz}\) is defined by \(\history_t (\tsample) \Let \hist_t = (\stt{}{0},\cont{}{0},\stt{}{1},\cont{}{1},\ldots,\stt{}{t})\).
\end{itemize}
Define a \emph{filtration} \(\bigl(\tsigalg_t \bigr)_{t \in \Nz}\) of \(\Borelsigalg{(\tss)}\) by \(\tsigalg_t \Let \sigma(\history_t)\) such that \(\tsigalg_t \subset \tsigalg_{t+1}\) for every \(t \in \Nz\). Observe that we have \[\Borelsigalg{(\tss)}= \bigvee_{t=0}^{+\infty} \tsigalg_t,\;\text{where }\tsigalg_t \Let \sigma(\history_t).\]
\textcolor{black}{Given an initial distribution \(\idist\) on \(\Borelsigalg{(\stsp)}\)} and an admissible policy \(\policy = (\policy_t)_{t \in \Nz}\), there exists \cite[Proposition C.10 and Remark C.11]{ref:HLJL-12} a probability measure \(\tRR^{\policy}_{\idist}\) on the measurable space \((\tss,\Borelsigalg{(\tss)})\) such that for each \(t \in \Nz\), and for all \(S_0\in\Borelsigalg(\stsp)\), \(T\in\Borelsigalg(\admact)\), and \(S_1\in\Borelsigalg(\stsp)\)
 \begin{equation}
     \label{eq:updf}
     \begin{aligned}
         \begin{cases}
              \tRR^\pi_\idist\probof[\big]{\proc_0\in S_0} = \idist(S_0), \\
		      \tRR^\pi_\idist\cprobof[\big]{\control_t\in T \given \history_t} = \policy_t\bigl( \history_t; T\bigr), \; \text{and}\\ 
		      \tRR^\pi_\idist\cprobof[\big]{\proc_{t+1}\in S_1\given \history_t, \control_t} = \tRR^\pi_\idist\cprobof[\big]{\proc_{t+1}\in S_1\given \proc_t, \control_t}.  
         \end{cases}
     \end{aligned}
 \end{equation}
 The stochastic process 
 \begin{equation}
 \label{eq:mdp_cps}
    \Bigl(\tss,\Borelsigalg{(\tss),\bigl(\tsigalg_t \bigr)_{t \in \Nz}},\tRR^{\policy}_{\idist}, (\proc_t)_{t \in \Nz} \Bigr)
 \end{equation}
 is called a \emph{Markov decision process.} 
}

\subsection{The concept of separability}
\label{subsec:concept_sep}
The germ of the  idea of separability traces back to the theory of disjoint dynamical systems \cite{ref:HF-67},\footnote{The idea of the disjointedness of the support of two measures is motivated from \cite{ref:LevPelPer-15}.} but in what follows we formally define the term separability for completeness in the context of security of CPSs. 
We restrict our attention to sample space \(\ps \Let \stsp^{\Nz}\) generated by \((\stt{}{t})_{t \in \Nz}\) with entries in \(\stsp\).

Let \(\Borelsigalg(\ps)\) denote the Borel \(\sigma \mbox{-}\)algebra generated on \(\ps\) and let \((\sigalg_t)_{t \in \Nz}\) denotes a filtration of \(\Borelsigalg(\ps)\) generated by the \emph{history process} \(\phistory \Let (\phistory_t)_{t \in \Nz}\) where \(\phistory_t (\sample) \Let \phist_t \Let \bigl(\stt{}{0},\stt{}{1},\ldots,\stt{}{t}\bigr)\). Denote by  \(\sample \in \ps\), a \emph{sample path} or a \emph{state trajectory}, where \(\sample \Let \bigl(\stt{}{0},\stt{}{1},
\stt{}{2},\ldots\bigr)\).

{\color{black}With this information, let us distinguish between two classes of sample paths:  
\begin{enumerate}
    \item \label{it:honest}\(\sample = (\stt{}{0},\stt{}{1},\ldots)\) that are generated under the influence of an \emph{honest policy} \(\policy \Let (\policy_t)_{t \in \Nz}\) --- an honest policy \((\policy_t)_{t \in \Nz}\) refers to the nominal or \emph{no-attack} condition
    where all the actuators are honest, and the distribution of these paths is \(\PP^{\policy}_{\idist}\);
    
    \item  \label{it:corrupt}\(\csample = (\cstt{}{0},\cstt{}{1},\ldots)\) that are generated under the influence of a \emph{corrupt} or \emph{attack policy} \(\cpolicy \Let (\cpolicy_t)_{t \in \Nz}\) --- a corrupt or attack policy \((\cpolicy_t)_{t \in \Nz}\) refers to the case of \emph{at least} one actuator being compromised, i.e., under attack, and the distribution of these paths is denoted by \(\PP^{\cpolicy}_{\idist}\).\footnote{The readers are referred to Appendix \ref{appen:eom} for the construction of the measures \(\PP^{\policy}_{\idist}\) and \(\PP^{\cpolicy}_{\idist}\).}
\end{enumerate}
Note that \(\cpolicy\) is \textbf{not} necessarily a complete replacement of the honest policy \(\policy\). Instead, the \emph{corrupt policy} \(\cpolicy\) depicts \emph{every arbitrarily clever attacks} in a general framework. Specifically, certain components of \(\cpolicy\) may act differently from that of \(\policy\) at some instances, depending on the attacker. We highlight that no specific assumptions on \(\cpolicy\) have been imposed at this stage.}

We seek to \emph{separate} the aforementioned classes of sample paths, i.e., to distinguish between the sample paths \(\sample\) and \(\csample\) in order to detect the presence of malicious actuators in the CPS \eqref{eq:mdp_cps}, to which end we introduce the following definition:
\begin{defn}
\label{def:separator}
A Borel measurable map \(\sep: \ps \lra \aset[\big]{0,1}\) is a \emph{separator} if
\begin{itemize}[label=\(\circ\)]
    \item for \(\PP^{\policy}_{\idist}\mbox{-}\)\emph{almost every} \(\sample\) we have \(\sep(\sample) = 0\), and 
    \item for \(\PP^{\cpolicy}_{\idist}\mbox{-}\)\emph{almost every} \(\csample\) we have \(\sep(\csample) = 1\).
\end{itemize}
\textcolor{black}{The cyber-physical system \eqref{eq:mdp_cps} \emph{admits separation} if there exists a separator \(\sep\) such that \(\sep(\sample) = 0\) and \(\sep(\csample) = 1\); therefore, the CPS \eqref{eq:mdp_cps} then has the \emph{separability} property. }
\end{defn}
\begin{rem}
\label{rem:inter_def}
\textcolor{black}{Consider the infinitely long sample paths \(\tsample_1\) and \(\tsample_2\), where \(\tsample_1\) corresponds to the \emph{no-attack} condition and \(\tsample_2\) corresponds to the \emph{attack} condition (when \emph{at least} one actuator is under attack), the above definition concerns with the existence of a \emph{Boolean classifier} \(\sep\) such that if \(\tsample_1\) and \(\tsample_2\) are fed into \(\sep\), we get \(0\) and \(1\), respectively. Note that the existence of such a separator \(\sep\) is independent of the specific (infinitely long) sample paths \(\tsample_1\) and \(\tsample_2\) fed to it. As a result and as mentioned in \S\ref{sec:intro}, in practice the above idea is not verifiable in an \emph{online} fashion: firstly, it requires infinitely long memory, and secondly, it addresses \emph{posterior} properties in the sense that the classification of the two infinitely long sample paths with different statistics can only be done \emph{after} they have been generated.}

\end{rem}
\textcolor{black}{Clearly, a separator in the sense of Definition \ref{def:separator} may not always exist. One of the natural questions to ask, therefore, is \textit{under what conditions can the existence/non-existence of a separator be asserted.} We address this important question in the context of the CPS \eqref{eq:mdp_cps} (and more specifically, in the context of stochastic linear time-invariant systems in \S\ref{sec:mrg} ahead) in this article and provide asymptotic guarantees for the existence/non-existence of a separator.}

{\color{black}Let us establish conditions for the detection of the presence of malicious components in \eqref{eq:mdp_cps} from the notion of separability in the sense of Definition \ref{def:separator}. To that end, we define:
 \begin{itemize}[label=\(\circ\)]
     \item  \(\safe \Let \aset[\big]{\sample \suchthat \sep(\sample) = 0}\) and
     \item  \(\unsafe \Let \aset[\big]{\sample \suchthat \sep(\sample) = 1}\).
 \end{itemize}
 Observe that the sets \(\safe\) and \(\unsafe\) denote the supports of the measures \(\PP^{\policy}_{\idist}\) and \(\PP^{\cpolicy}_{\idist}\), respectively.
 \begin{prop}
     \label{prop:det_from_sep}
     Consider the CPS \eqref{eq:mdp_cps} along with its associated data and notations. Assume that a separator exists. If \(\unsafe \neq \emptyset\), then there exists at least one malicious component in \eqref{eq:mdp_cps}.
\end{prop}
\begin{pf}
    Observe that if \(\unsafe \neq \emptyset\), then there exists an element \(\ol{\tsample}_s \in \unsafe\) such that \(\sep(\ol{\tsample}_s) = 1\). Then from Definition \ref{def:separator} and the fact that \(\sep(\ol{\tsample}_s) = 1\) we conclude that the sample path \(\ol{\tsample}_s\) is generated under the influence of a corrupt policy \(\cpolicy\), which immediately asserts the presence of at least one malicious component in the CPS \eqref{eq:mdp_cps}. \hfill \(\blacksquare\)
\end{pf} 

The case \(\unsafe = \emptyset\) refers to the existence of a separator \(\sep\) wherein the sets \(\ps\) and \(\emptyset\) are considered to be separated despite not being accounted in Definition \ref{def:separator}. In our context, Definition \ref{def:separator} only concerns the existence of a \emph{non-trivial} separator (in the sense that \(\unsafe \neq \emptyset\)) and will be termed as a \emph{separator} in the sequel. In view of Proposition \ref{prop:det_from_sep} the problem of detecting the presence of malicious components in the CPS then boils down to establishing sufficient conditions for the existence and non-existence of a separator, which is the main topic of \S\ref{subsec:chief_obs}.
}

\subsection{Chief observation}
\label{subsec:chief_obs}
Here we state the first chief observation of this article. It provides a set of necessary and sufficient conditions for the existence of a separator. Let us define the measures \(\QQ_n\) and \(\cQQ_n\) on the measurable space \((\ps,\sigalg_n)\): for any \(E \in \sigalg_n\) of the form \(E = S_0\times S_1\times \cdots\times S_{n-1} \times S_n\) with Borel sets \(S_i \in \Borelsigalg(\stsp)\) for each \(i \in \Nz\), we have
\begin{align}
    \label{e:marginal}
        \QQ_n \probof[\big]{E} &= \PP^{\policy}_{\idist}(\proc_o \in S_0,\proc_1 \in S_1, \ldots,\proc_n \in S_n), 
\end{align}
and 
\begin{align}
    \label{e:cmarginal}
        \cQQ_n \probof[\big]{E} &= \PP^{\cpolicy}_{\idist}(\proc_o \in S_0,\proc_1 \in S_1, \ldots,\proc_n \in S_n).
\end{align}
Note that the measures \(\QQ_n\) and \(\cQQ_n\) are the restrictions of \(\PP^{\policy}_{\idist}\) and \(\PP^{\cpolicy}_{\idist}\), respectively, to \(\sigalg_n\). See Appendix \ref{appen:eom} for an assertion on the existence of the measures \(\QQ_n\) (and \(\cQQ_n\)) and \(\PP^{\policy}_{\idist}\) (and \(\PP^{\cpolicy}_{\idist}\)).
\begin{thm}
	\label{p:key}
	\textcolor{black}{Consider the CPS \eqref{eq:mdp_cps} along with its associated data and notations. Let \(\idist\) denote an initial probability distribution on \(\stsp\), and let \(\policy \Let (\policy_t)_{t\in\Nz}\) and \(\cpolicy \Let (\cpolicy_t)_{t\in\Nz}\) denote the control policies corresponding to the case of all honest actuators and to the case of at least one malicious actuator, respectively. Let \(\sample\) and \(\csample\) be the sample paths generated under the influence of \(\policy\) and \(\cpolicy\), respectively. Assume that \(\QQ_n \abscont \cQQ_n \) for each \(n \in \Nz\). Then there exists a separator in the sense of Definition \eqref{def:separator} if and only if  }
	{\color{black}\begin{equation}\label{e:key}
		\frac{\odif{\QQ_n}}{\odif{\cQQ_n}} \xrightarrow[n\to+\infty]{\PP^{\cpolicy}_{\idist}-\text{a.s.}} 0,
	\end{equation}
    where \(\frac{\odif{\QQ_n}}{\odif{\cQQ_n}}\) represents the Radon-Nikodym derivative of \(\QQ_n\) relative to \(\cQQ_n\).}
\end{thm}
{\color{black}
\begin{pf}
    It suffices to establish that mutual singularity of the probability measures \(\PP^{\policy}_{\idist}\) and \(\PP^{\cpolicy}_{\idist}\) is equivalent to the existence of a separator, defined in Definition \ref{def:separator}.\footnote{\textcolor{black}{Refer to Appendix \ref{appen:prob basics} for the definition of mutual singularity of two measures.}} 
    
    Assume that \eqref{e:key} holds. Invoking Jessen's theorem \cite[Theorem 5.2.20]{ref:Str-11} shows that the limiting probability measures \(\PP^{\policy}_{\idist}\) and \(\PP^{\cpolicy}_{\idist}\) are mutually singular, i.e., \(\PP^{\policy}_{\idist} \perp \PP^{\cpolicy}_{\idist}\), which implies that there exist a Borel measurable set \(S \in \Borelsigalg{(\ps)}\) such that \[ \PP^{\policy}_{\idist}\Bigl(\bigl(\proc_0, \proc_1,\ldots, \proc_n,\dots\bigr) \in S\Bigr)=1\] and \[ \PP^{\cpolicy}_{\idist}\Bigl(\bigl(\proc_0, \proc_1,\ldots, \proc_n,\ldots\bigr) \in S^{\mathsf{c}} \Bigr)=1.\] Observe that the random variable \[\ps\ni \sample \mapsto \indic{\aset[]{(\proc_0,\proc_1,\ldots)\in S^{\mathsf{c}}}}(\sample)\] 
is a Borel measurable function and by definition, it separates the support of the probability measure \(\PP^{\policy}_{\idist}\) and \(\PP^{\cpolicy}_{\idist}\) asserting the existence of a separator. 

We now establish the converse statement. Recall the definitions of the sets \(\safe\) and \(\unsafe\):
 \begin{itemize}[label=\(\circ\)]
     \item  \(\safe \Let \aset[\big]{\sample \suchthat \sep(\sample) = 0}\) and
     \item  \(\unsafe \Let \aset[\big]{\sample \suchthat \sep(\sample) = 1}\).
 \end{itemize}

If there exists a Borel measurable function \(\sep: \ps \lra \aset[\big]{0,1}\) such that for \(\PP^{\policy}_{\idist}\mbox{-}\)\emph{almost every} \(\sample\) we have \(\sep(\sample) = 0\) and for \(\PP^{\cpolicy}_{\idist}\mbox{-}\)\emph{almost every} \(\csample\) we have \(\sep(\csample) = 1\). We observe that \(\ps = \safe \cup \unsafe\), and the sets \(\safe\) and \(\unsafe\) are disjoint, i.e., \(\safe \cap \unsafe = \emptyset\). From Definition \ref{def:separator}, it follows that \(\PP^{\cpolicy}_{\idist} (\unsafe) = 1\) and \(\PP^{\policy}_{\idist} (\safe) = 1\), asserting that \(\PP^{\policy}_{\idist}\) and \(\PP^{\cpolicy}_{\idist}\) are mutually singular. Invoking Jessen's theorem \cite[Theorem 5.2.20]{ref:Str-11} immediately gives us 
\[
    \frac{\odif{\QQ_n}}{\odif{\cQQ_n}} \xrightarrow[n\to+\infty]{\PP^{\cpolicy}_{\idist}-\text{a.s.}} 0,
\]
which concludes the proof.    \hfill \(\blacksquare\)
\end{pf}
}
{\color{black}
\begin{rem}
    \label{rem:key remark}
    Theorem \ref{p:key} is a foundational result and it sets the bedrock for a general framework for \textbf{arbitrarily clever actuator attack detection} of the CPS \eqref{eq:mdp_cps} under the setting of Markov decision process. As mentioned before, we do not impose any restrictions and limitations on the \emph{attack} or \emph{corrupt policy} \(\cpolicy\), which further contributes to its vast generality.     
\end{rem}

Theorem \ref{p:key} furnishes us with the mathematical tools to establish negative results that point toward the fundamental limitations and trade-offs that one must encounter while designing secured control systems. We study this aspect in detail in \S\ref{sec:mrg} for the case of stochastic linear time-invariant systems with Gaussian noise.
}

\section{Sufficient conditions for the (non)existence of a separator}
\label{sec:mrg}
\textcolor{black}{Here we restrict our attention to stochastic linear cyber-physical systems (CPSs) and obtain crisp sufficient conditions for the existence and non-existence of a separator in the sense of Definition \ref{def:separator}. We stipulate that the  \emph{honest policy} \(\policy\) and the \emph{corrupt policy} \(\cpolicy\) belong to class of randomized policies, denoted by \(\pclass{R}\). The defense mechanism adopted by the CPS involves injecting a random perturbation, termed as \emph{private excitation}, with suitable statistics, on top of the ``nominal'' control actions. In this section we investigate the case of stochastic linear CPSs with Gaussian private excitation and process noise.
}
\subsection{System description}
\label{sec:ls}
Consider a cyber-physical system (CPS) with \(\node\) nodes where the \(i^{\mathrm{th}}\) node represents agent \(i\) described by a scalar stochastic linear  system 
\begin{equation}
\label{eq:lsyss}
    \stt{i}{t+1} = a_{ii} \stt{i}{t} + \sum_{j \in \ngb_i}a_{ij} \stt{j}{t} + b_{i}\cont{i}{t} + \pnoise{i}{t} \qquad \text{where\,\,} i=1,2,\ldots,\node,
\end{equation}
with the following data: for $i,j = 1,2,\ldots,\node$,
\begin{enumerate}[label=\textup{(\ref{eq:lsyss}-\alph*)}, leftmargin=*, widest=b, align=left]
    \item \label{it:state} \(\stt{i}{t} \in \Rbb\) is the state of the \(i^{\mathrm{th}}\) agent at time \(t\);
    
    \item \label{it:control} \(\cont{i}{t} \in \Rbb\) denotes the control action applied to the \(i^{\mathrm{th}}\) agent at time \(t\);
    
    \item \label{it:sp} the scalars \(a_{ij} \in \Rbb\) and \(b_{i} \in \Rbb\) denote the system and control parameters corresponding to the \(i^{\text{th}}\) agent of the CPS; 

    \item \label{it:inneighbour}  the set \(\ngb_i\) denotes the set of in-neighbours of the \(i^{\mathsf{th}}\) agent in the network;\footnote{In-neighbors of the \(i^{\mathsf{th}}\) agent refers to the set of all neighboring agents that can directly influence the evolution of the \(i^{\mathsf{th}}\) agent.}
    
    \item \label{it:pn} \((\pnoise{i}{t})_{t \geqslant 0}\) is a sequence of independent and identically distributed (i.i.d) random variables with some known distribution, denoting the process noise corresponding to the \(i^{\mathrm{th}}\) agent;
    
    \item \label{it:ic} \(\stt{i}{0}\) is the initial state of the \(i^{\mathsf{th}}\) agent with some known distribution \(\idist^i\) and is assumed to be independent of the sequence \((\pnoise{i}{t})_{t \geqslant 0}.\)
\end{enumerate}
\textcolor{black}{This model of CPS is motivated from \cite{ref:MAR-AGA-13} where each agent in the network is modeled as an integrator.} 

The overall system of \(\node\) agents can be expressed in the following compact form:
 \begin{equation}
     \label{eq:lsyscompactt}
     \stt{}{t+1}=A\stt{}{t}+B\cont{}{t}+\pnoise{}{t},
 \end{equation}
 where \(A= [a_{ij}] \in \Rbb^{\node \times \node}\), \(B=\mathsf{diag} (b_1, b_2, \ldots,b_{\node}) \in \Rbb^{\node \times \node}\) denote the state matrix and the input matrix, respectively.\footnote{The current setting assumes that every agent is equipped with a local controller. In many applications it is possible that an agent, e.g., the \(i^{\text{th}}\) agent, does not possess a local controller which will result in \(b_i = 0\); those cases are within the scope of this article.} \textcolor{black}{The process noise \((\pnoise{}{t})_{t \in \Nz}\) where \(\pnoise{}{t}=\bigl(\pnoise{1}{t} \;\pnoise{2}{t}\; \cdots\; \pnoise{\node}{t}\bigr)^{\top}\) for \(t \in \Nz\), is a Gaussian random process with i.i.d. entries and \(\pnoise{}{0} \sim \grv(0,\var_{\pnoise{}{}})\) with \(\var_{\pnoise{}{}} = \var_{\pnoise{}{}}^{\top} \succ 0\).}

 Recall that the admissible state-space and the admissible action space are denoted by \(\stsp \subset  \Rbb^{\node}\) and \(\admact\subset \Rbb^{\node}\), respectively.
The control action \(\cont{}{t}\) at time \(t\) in the presence of all honest agents is chosen according to a pre-determined admissible control strategy or policy \(\policy\). \textit{Here we focus on the class of \emph{randomized} policies, both history dependent and Markovian.}



\begin{lem}
   \label{lem:wrpm}
   Consider the CPS \eqref{eq:lsyss}/\eqref{eq:lsyscompactt} along with its associated data \ref{it:state}-\ref{it:ic}. Let \(\idist\) be an arbitrary initial distribution on \(\stsp\), i.e., \(\PP^{\policy}_{\idist} \probof[\big]{\proc_0 \in S_0} = \idist \probof[\big]{S_0}\), and let \((\stt{}{n})_{n \in \Nz}\) denote a trajectory generated by the system \eqref{eq:lsyss}/\eqref{eq:lsyscompactt} under \(\policy \in \pclass{}\). Let \(E \in \sigalg_n\) be a Borel set of the form  \(E = S_0\times S_1 \times\cdots\times S_{n-1} \times S_n\) where \(S_i \in \Borelsigalg(\stsp)\) for every \(i=1,\ldots,n\). Then, for each \(n \in \Nz\), the measure defined in \eqref{e:marginal} takes the form 
   \begin{equation}
       \label{eq:lem:wrpm}
       \begin{aligned}
          \QQ_n &(E) = \int_{S_0} \idist(\odif{y_o})\int_{S_1} \PP^{\policy}_{\idist} \cprobof[\big]{\proc_1 \in \odif{y_1} \given \phist_0} \cdots \,\\& \cdots \int_{S_{n-1}} \PP^{\policy}_{\idist} \cprobof[\big]{\proc_{n-1} \in \odif{y_{n-1}} \given \phist_{n-2}}\int_{S_n} \PP^{\policy}_{\idist} \cprobof[\big]{\proc_n \in \odif{y_n} \given \phist_{n-1}},
       \end{aligned}
   \end{equation}
   where \(\phist_n \Let \bigl(y_0,\ldots,y_{n}\bigr)\) denotes the entire history of the CPS until the \(n^{\text{th}}\) time instant.
\end{lem}
 Moreover, we write the measure \eqref{eq:lem:wrpm} as
   \[
   \QQ_n\probof[\big]{\odif{y_0},\ldots, \odif{y_n}} =  \idist(\odif{y_o}) \prod_{t=1}^n \PP^{\policy}_{\idist} \cprobof[\big]{\proc_t \in \odif{y_t} \given \phist_{t-1}}.  
   \]
A proof of Lemma \ref{lem:wrpm} is given in Appendix \ref{appen:proof_auxx}. An immediate corollary of Lemma \ref{lem:wrpm} is:
\begin{cor}
\label{c:wrpm}
Suppose the hypotheses of Lemma \ref{lem:wrpm} hold. If the policy \(\policy\) belongs to the class of randomized-Markovian policies \(\pclass{RM}\), then for every \(n \in \Nz\)
\begin{align}\label{cor:wrpm}
    \QQ_n & (E) = \int_{S_0} \idist(\odif{y_o})\int_{S_1} \PP^{\policy}_{\idist} \cprobof[\big]{\proc_1 \in \odif{y_1} \given y_{0}} \cdots \,\nn\\& \cdots \int_{S_{n-1}} \PP^{\policy}_{\idist} \cprobof[\big]{\proc_{n-1} \in \odif{y_{n-1}} \given y_{n-2}}\int_{S_n} \PP^{\policy}_{\idist} \cprobof[\big]{\proc_n \in \odif{y_n} \given y_{n-1}}.
\end{align}
\end{cor}
\textcolor{black}{The proof of Corollary \ref{c:wrpm} follows directly from Lemma \ref{lem:wrpm} by substituting \((\phist_n)_{n \in \Nz}\) with \((\stt{}{n})_{n \in \Nz}\) under the assumption that only the current history is taken into account.}

\subsection{Attack description}
\label{subsec:types}
{\color{black}
Theorem \ref{p:key} naturally include \emph{arbitrarily clever actuator attacks}. In our setting malicious actuators can collude with each other, share critical information and devise strategies to tamper with the honest actions of the CPS, observe the states \emph{silenly} and learn the system parameters if the underlying dynamics is hidden, and at the same time remain stealthy. It encompasses a large class of randomized attacks that are \emph{history dependent} and, in particular, \emph{Markovian}, i.e., \(\policy, \cpolicy \in \pclass{RH}\) or \(\policy, \cpolicy \in \pclass{RM}\); which is the main topic of \S\ref{subsec:policy classess}. 

Moreover, we draw attention to the fact that the precise functional forms of the policies \(\policy\) and \(\cpolicy\), as long as they belong to the class of all randomized policies, are inconsequential in what follows. In other words, the results obtained in \S\ref{subsec:policy classess} are agnostic to the specific functional nature of the randomized policies played by the actuators (honest and corrupt). We direct the readers to Fig. \ref{fig:attack model} where the attack-prone components are marked in red.
\begin{figure}[!htpb]
\centering
\includegraphics[scale = 0.23]{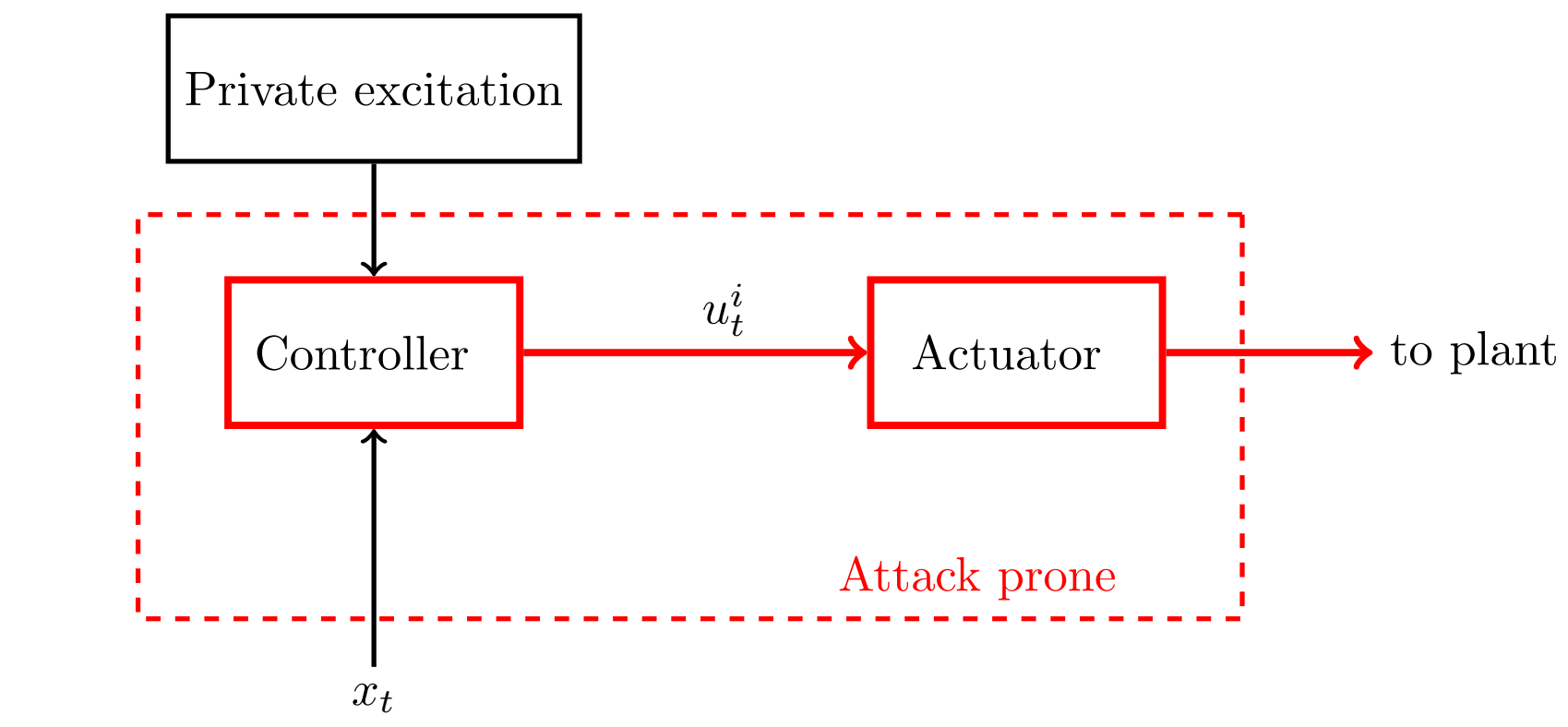}
\caption{Diagram of the attack prone components corresponding to agent \(i\).}
\label{fig:attack model}
\end{figure}

\begin{rem}
\label{rem:import}
In this article a malicious attacker can only hijack a subset of the actuators. We do not consider the case of moving hijackers that change their locations. If all the actuators in the CPS are compromised, then the attacker will have access to every control input and the entire history of the system, and with complete information on the system model, it can remain stealthy forever. Such an omnipresent attacker is not considered in this work since no detection methods can be leveraged to identify its presence in such a case. This is also a point of departure from arbitrary sensor attacks where even if all the sensors are compromised, the actuators remain secured and can take appropriate measures.
\end{rem}

Specifically, an attacker may possess one or more of the following capabilities:
\begin{itemize}[label = \(\circ\)]
    \item It could completely overhaul the nominal control action by directly tampering with it.
    \item It could potentially hijack those communication channels that are highlighted in red in Fig. \ref{fig:attack model}.
    \item It could launch physical attacks on the actuators and their subsystems.\footnote{Physical attacks on the actuator refer to those cases where the attackers can permanently or temporarily damage the actuators or any associated subsystems. }
\end{itemize}
Note that all these capabilities of the attacker are jointly represented by the corrupt policy \(\cpolicy\) and the corrupt control action \(\ccont{}{t}\).\footnote{Some  specific instances include cutting off the communication channel between the controllers and the actuators or adding an exogenous input signal, popularly known as ``DoS attacks'' and ``false data injection attacks''.} While silently observing the output channel is allowed, we reiterate that tampering or modifying the output communication channel is not considered in this work.

To identify such malicious attacks, we adopt the approach of injecting private excitations by the honest components with suitable statistics, superimposed on the pre-determined control actions. It admits the realization
\begin{equation}
    \label{eq:nom_control}
    \cont{}{t} = \action_t + \pe{}{t} \quad \text{for every }t \in \Nz,
\end{equation}
where
\begin{itemize}[label = $\circ$]
    \item \(\action_t\) corresponds to the pre-determined control action at time \(t\) designed under normal circumstances (no-attack situation), and
    \item \(\pe{}{t}\) denotes the private excitation injected into the system at time \(t\), to add a layer of security against malicious attacks.
\end{itemize}
}
\begin{assum}\label{assum:gaussian case}
The following assumptions are imposed on the private excitation throughout the sequel:
\begin{enumerate}[label=\textup{(\ref{assum:gaussian case}-\alph*)}, leftmargin=*, widest=b, align=left]
    \item \label{it:gaussian case-1} The private excitation denoted by \((\pe{}{t})_{t \in \Nz}\) is a Gaussian random process with i.i.d. entries and \(\pe{}{0} \sim \grv(0,\var_{\pe{}{}})\); its distribution is made public and is \emph{known} to the malicious attackers. Moreover, it is assumed to be independent of all other random vectors in the system.  
    \item  The components of \(\pe{}{t}\) are uncorrelated for every \(t\), i.e., for every distinct \(i,j \in \aset[]{1,\ldots,\node}\) we have \(\pe{i}{t}\) uncorrelated with \(\pe{j}{t}\); that is, \(\var_{\pe{}{}}\) is a diagonal matrix.
\end{enumerate}
\end{assum}


\subsection{Sufficient conditions for (non)separability under \(\pclass{RH}\)}
\label{subsec:policy classess}
The set \(\pclass{RH}\) consists of \emph{randomized} and \emph{history dependent policies}. Recall from Definition \ref{def:policy} in Appendix \ref{appen:prob basics} that such a policy is a sequence \((\policy_t)_{t \in \Nz}\) of stochastic kernels (Definition \ref{def:sk} in Appendix \ref{appen:prob basics}) on the set \(\admact\) given the entire past \(\hist_t\) satisfying \(\policy_t\bigl(\hist_t;\admact\bigr) = 1.\)

Let \(\phist_t \Let (\stt{}{0},\ldots,\stt{}{t})\) and \(\cphist_t \Let (\cstt{}{0},\ldots,\cstt{}{t})\) denotes the history until the \(t^{\text{th}}\) time instant, under \(\policy\) and \(\cpolicy\), respectively. Then the entries of \((\action_t)_{t \in \Nz}\) in \eqref{eq:nom_control} take the form
\(\action_t = \pmap{t}(\phist_t)\), where each policy map \(\pmap{t}\) is a Borel measurable function. 
\textcolor{black}{Similarly, the entries of \((\caction_t)_{t \in \Nz}\) under the \emph{attack situation} take the form \(\caction_t = \cpmap{t}(\cphist_t)\), where \(\cpmap{t}\) is a Borel measurable function. Note that, as mentioned in \S\ref{subsec:types}, the specific functional forms of \(\pmap{t}(\cdot)\) and \(\cpmap{t}(\cdot)\) are not pertinent, and the results derived here are agnostic to them. The control input admitted by the agent \(i\) at time \(t\) is given by
\begin{equation}
\label{controlnoattack}
\cont{i}{t}=\pmap{t}^i(\phist_t)+ \pe{i}{t},
\end{equation}
where \(\pmap{t}^i(\cdot)\) corresponds to the \(i^{\text{th}}\) component of \(\pmap{t}(\cdot)\). Similarly, when the CPS is under attack, the input admitted by a malicious agent \(j\) at time \(t\) is written as
\begin{equation}
\label{controlattack}
  \ccont{i}{t}=\cpmap{t}^i(\cphist_t),
\end{equation}
where \(\cpmap{t}^i(\cdot)\) is the \(i^{\text{th}}\) component of \(\cpmap{t}(\cdot)\). }

Note that a malicious attacker can hijack a subset of the actuators, as mentioned in Remark \ref{rem:import}. Let \(M\) be the number of actuators hijacked by a malicious attacker. \textcolor{black}{For the sake of analysis, we stack the control actions generated by the malicious agents and the control actions generated by the honest agents separately.} That is, the control action  \(\ccont{}{t}\) generated by \(\cpolicy\) is partitioned as \(\bigl(\ccont{}{t,1}\;\cont{}{t,2}\bigr)^{\top}\), where \(\ccont{}{t,1}=\bigl(\ccont{1}{t}\; \ccont{2}{t}\; \cdots\; \ccont{M}{t}\bigr)^{\top}\) and \(\cont{}{t,2}=\bigl(\cont{M+1}{t}\; \cont{M+2}{t}\; \cdots\; \cont{N}{t}\bigr)^{\top}\), respectively, and the control action \(\cont{}{t}\) generated by \(\policy\) is partitioned as \(\big(\cont{}{t,1}\;\cont{}{t,2}\bigr)^{\top}\), where \(\cont{}{t,1}=\bigl(\cont{1}{t}\; \cont{2}{t}\; \cdots\; \cont{M}{t}\bigr)^{\top}\). 

As a consequence of \eqref{controlnoattack} and \eqref{controlattack}, it follows that the control actions generated by \(\policy\) and \(\cpolicy\), when the system \eqref{eq:lsyscompactt} operates under \emph{no-attack} condition and when it is under attack by the adversaries, respectively, are given by 
\begin{align}
    \label{eq:GRM policy}
    &\cont{}{t} = \begin{pmatrix} \cont{}{t,1}\\\cont{}{t,2}\end{pmatrix}=\begin{pmatrix}\pmap{t,1}(\phist_t) + \pe{}{t,1}\\\pmap{t,2}(\phist_t) + \pe{}{t,2}
    \end{pmatrix} \quad \text{and} \nn\\ &\ccont{}{t}=\begin{pmatrix} \ccont{}{t,1}\\\cont{}{t,2}\end{pmatrix} = \begin{pmatrix}\cpmap{t,1}(\cphist_t)\\\pmap{t,2}(\cphist_t) + \pe{}{t,2}
    \end{pmatrix}; 
\end{align}
here the policy maps are partitioned as 
\(\pmap{t,1}=\bigl(\pmap{t}^1\;\cdots\; \pmap{t}^M\bigr)^{\top}\), \(\pmap{t,2}=\bigl(\pmap{t}^{M+1}\;\cdots\; \pmap{t}^N\bigr)^{\top}\) and \(\cpmap{t,1}=\bigl(\cpmap{t}^1\;\cdots\; \cpmap{t}^M\bigr)^{\top}\), and the private excitation injected into the system is also partitioned as \(\pe{}{t,1}=\bigl(\pe{1}{t}\; \pe{2}{t} \;\cdots\; \pe{M}{t}\bigr)^{\top}\) and \(\pe{}{t,2}=\bigl(\pe{M+1}{t}\; \pe{M+2}{t} \;\cdots\; \pe{N}{t}\bigr)^{\top}\).  
{\color{black}\begin{assum}
    \label{assum:stabilizability}
    We stipulate that each agent is influenced directly or indirectly by at least one of the honest actuators.
\end{assum}
 Assumption \ref{assum:stabilizability} refers to the necessary condition imposed on the underlying network of the CPS \eqref{eq:lsyss}/\eqref{eq:lsyscompactt} to guarantee that the private excitation injected by the honest agents reach the attacked actuators, either directly or indirectly, so that \(\QQ_n \abscont \cQQ_n\) for each \(n \in \Nz\), is ascertained. For more detail, we refer the readers to Example \ref{exam:graph2} ahead.
 In particular assumption \ref{assum:stabilizability} ensures that the CPS is at least stabilizable with respect to the honest actuators. 
}

Next, we provide sufficient conditions for the (non)existence of a separator in the sense of Definition \ref{def:separator} under randomized and history dependent policies:
{\color{black}\begin{thm}
\label{prop:main}
Consider the CPS \eqref{eq:lsyss}/\eqref{eq:lsyscompactt} along with its associated data \ref{it:state}-\ref{it:ic} and notations. Let \(M\) denotes the number of actuators in the network that are under attack when the CPS is hijacked by an adversary. Suppose that Assumptions \ref{assum:gaussian case} and \ref{assum:stabilizability} hold. Let \((\policy_t)_{t \in \Nz}\) and \((\cpolicy_t)_{t \in \Nz}\) denote the history dependent control policies corresponding to the case of \emph{all} honest actuators and to the case of attack, respectively. Define the quantity
\[\rate_n \Let \frac{\sum_{t=1}^n \eig{-1}{\min,t} \norm{\stt{}{t}-\meantr{\policy}{t}}^2}{\sum_{t=1}^n \ceig{-1}{\max,t} \norm{\cstt{}{t}-\meantr{\cpolicy}{t}}^2} \qquad \text{for }n \in \Nz, \]
where for each \(t \in \Nz\),
\begin{itemize}[leftmargin =*, label=\(\circ\)]
\item \(\meantr{\policy}{t}\) and \(\meantr{\cpolicy}{t}\) denote the conditional expectation of \(\stt{}{t}\) and \(\cstt{}{t}\)  given the entire history \(\phist_{t-1}\) and \(\cphist_{t-1}\), respectively;
\item \(\eig{}{\min,t}\) and \(\ceig{}{\max,t}\) denote the minimum and maximum eigenvalues corresponding to the conditional variance of \(\stt{}{t}\) and \(\cstt{}{t}\), given their entire history \(\phist_{t-1}\)  and \(\cphist_{t-1}\), respectively.  
\end{itemize}
Then we have the following assertions: 
\begin{enumerate}[label=\textup{(T\ref{prop:main}-\alph*)}, leftmargin=*, widest=b, align=left]
    \item \label{it:fa2} If \(\rate_n \xrightarrow[n\to+\infty]{\PP^{\cpolicy}_{\idist}-\text{a.s.}} 1\), then the CPS \eqref{eq:lsyss}/\eqref{eq:lsyscompactt} does not admit separation in the sense of Definition \ref{def:separator}.
    
    \item \label{it:faaaa2} If there exist a constant \(M>0\) such that 
    \begin{equation*}
        \begin{aligned}
            \begin{cases}\sup_{n \in \Nz}\rate_n \leqslant M &\PP^{\cpolicy}_{\idist}-\text{a.s.},\\
            \limsup_n \sum_{t=1}^n \ceig{-1}{\max,t} \norm{\cstt{}{t}-\meantr{\cpolicy}{t}}^2 < +\infty & \PP^{\cpolicy}_{\idist}-\text{a.e.}\mbox{-}\sample, \end{cases}
        \end{aligned} 
    \end{equation*}
    then the CPS \eqref{eq:lsyss}/\eqref{eq:lsyscompactt} does not admit separation in the sense of Definition \ref{def:separator}.

    \item \label{it:ta2} If the following hypotheses hold:
    \begin{equation*}
        \begin{aligned}
            \begin{cases}\rate_n \xrightarrow[n\to+\infty]{\PP^{\cpolicy}_{\idist}\mbox{-}\text{a.s.}} + \infty, \, \text{and}\\
            {\sum_{t=1}^n \ceig{-1}{\max,t} \norm{\cstt{}{t}-\meantr{\cpolicy}{t}}^2 \xrightarrow[]{n\to+\infty} +\infty}\quad \PP^{\cpolicy}_{\idist}-\text{a.e.}\mbox{-}\sample,\end{cases}
        \end{aligned}
    \end{equation*} 
    then there exists a separator \(\sep\) and the CPS \eqref{eq:lsyss}/\eqref{eq:lsyscompactt} admits separation in the sense of Definition \ref{def:separator}.
\end{enumerate}
\end{thm}

The proof of Theorem \ref{prop:hmain} is based on the fact that given the history \(\phist_{t-1}\), the random vector \(\proc_t\) is conditionally Gaussian under Assumption \ref{assum:gaussian case}, which is the main thesis of next Lemma.
\begin{lem}
  \label{lem:rvg}
  Consider the CPS \eqref{eq:lsyss}/\eqref{eq:lsyscompactt} with its associated data \ref{it:state}-\ref{it:ic} and notations, and assume that the hypotheses of Theorem \ref{prop:main} hold. Then for every \(t \in \Nz\), \(\proc_t\) is a conditionally Gaussian random vector under \(\policy_{t-1}\) and \(\cpolicy_{t-1}\) given \(\phist_{t-1}\) and \(\cphist_{t-1}\), respectively.
\end{lem}
Proofs of Lemma \ref{lem:rvg} and Theorem \ref{prop:main} have been relegated to Appendix \ref{appen:proof_auxx} in \S\ref{appen:main_proof}.}

\subsection{Sufficient conditions for (non)separability under \(\pclass{RM}\)}
\label{subsubsec:rh control}
 {\color{black}The set \(\pclass{RM}\) consists of \emph{randomized} and \emph{Markovian} policies --- such a policy is a sequence \((\policy_t)_{t \in \Nz}\) of stochastic kernels on the set \(\admact\) given \(\stt{}{t}\) satisfying \(\policy_t(\stt{}{t};\admact) = 1.\) Here we adhere to various notations adopted in \S\ref{subsec:policy classess}. 
The next result is an immediate corollary of Theorem \ref{prop:main} that provides sufficient conditions for the (non)existence of a separator when \(\pclass{RM}\) is considered in place of \(\pclass{RH}\).
\begin{cor}
\label{prop:hmain}
Suppose that the hypotheses of Theorem \ref{prop:main} hold. Let all the actuators employ \emph{randomized Markovian policies}, i.e., \(\policy, \cpolicy \in \pclass{RM}\). Define the quantity
\begin{equation}\label{eq:rate}
    \rate_n \Let \frac{\sum_{t=1}^n \eig{-1}{\min,t} \norm{\stt{}{t}-\meantr{\policy}{t}}^2}{\sum_{t=1}^n \ceig{-1}{\max,t} \norm{\cstt{}{t}-\meantr{\cpolicy}{t}}^2} \qquad \text{for }n \in \Nz,
\end{equation}
where  for each \(t \in \Nz,\) 
\begin{itemize}[leftmargin=*, label=\(\circ\)]
\item \(\meantr{\policy}{t}\) and \(\meantr{\cpolicy}{t}\) denote the conditional expectation of \(\stt{}{t}\) and \(\cstt{}{t}\) given the current history \(\stt{}{t-1}\) and \(\cstt{}{t-1}\), respectively, and
\item \(\eig{}{\min,t}\) and \(\ceig{}{\max,t}\) denote the minimum and maximum eigenvalues corresponding to the conditional variance of \(\stt{}{t}\) and \(\cstt{}{t}\) given the current history \(\stt{}{t-1}\) and \(\cstt{}{t-1}\), respectively. 
\end{itemize} Then, the following assertions hold:
\begin{enumerate}[label=\textup{(C\ref{prop:hmain}-\alph*)}, leftmargin=*, widest=b, align=left]
    \item \label{it:cor-markov} The closed loop process \(\proc = (\proc_t)_{t \in \Nz}\) \textcolor{black}{under randomized Markovian policies} is Markovian.

    \item \label{it:cor-gaussian} for every \(t \in \Nz\), \(\proc_t: \ps \lra \Rbb^{\node}\) is a conditionally Gaussian random vector under the policies \((\policy_t)_{t \in \Nz}\) and \((\cpolicy_t)_{t \in \Nz}\) given \(\stt{}{t-1}\) and \(\cstt{}{t-1}\), respectively.

    \item \label{it:cor-th2} The assertions \ref{it:fa2}-\ref{it:ta2} in Theorem \ref{prop:main} hold.
\end{enumerate}
\end{cor}
Please see \S\ref{appen:cor_proof} in Appendix \ref{appen:proof_auxx} for a proof of Corollary \ref{prop:hmain}.
}
\subsection{Examples}
\label{subsec:examples}
{\color{black}Let us understand the effect of the network structure and the choice of suitable private excitations with the help of the following simple examples, which further justify the assumptions made in Theorem \ref{prop:main} and Corollary \ref{prop:hmain}.
\begin{example}[Effect of the network structure]\label{exam:graph2}
For illustration, we consider a CPS consisting of \(\node = 2\) agents where the system matrix and input matrix are denoted by \(A = \diag(a_{11},a_{22})\) and \(B = \diag(b_1,b_2)\), respectively.\footnote{The discussion in Example \ref{exam:graph2} can be generalized to arbitrary \(\node\).}
Here \(a_{11}, a_{22} \neq 0\) and \(b_1,b_2 \neq 0\). We assume that:
  \begin{itemize}[leftmargin=*]
      \item \label{it:1} \(\mnode = 1,\) i.e., one of the actuators is malicious;
      \item \label{it:2}\((\pnoise{}{t})_{t\in \Nz}\) is a sequence of i.i.d. random vectors. Moreover, \(\bigl(\pnoise{1}{t}\bigr)_{t \in \Nz}\sim \mathsf{Unif}([-c,c])\) and \(\bigl(\pnoise{2}{t}\bigr)_{t \in \Nz} \sim \mathsf{Unif}([-c,c])\) for some \(c>0\), are independent of each other; 
      \item \label{it:3} the honest actuators employ state-feedback, i.e., the control actions under the no-attack situation and under the attack situation, are given by \(\cont{}{t} = \bigl(k_{t,1}\stt{1}{t} + \pe{}{t,1}, k_{t,2}\stt{2}{t} + \pe{}{t,2}\bigr)\) and  \(\ccont{}{t} = \bigl(\ccont{}{t,1},\, k_{t,2}\cstt{2}{t} + \pe{}{t,2}\bigr)\), respectively. Here \((\ccont{}{t})_{t \in \Nz} \sim \cpolicy \subset  \pclass{RM}\) denotes the attack scheme employed by the malicious actuator.
  \end{itemize}
Note that in this example, the overall system \eqref{eq:lsyscompactt} is decoupled implying that the agents do not communicate with each other. We fix the time \(n =1\). From \eqref{e:marginal} and \eqref{e:cmarginal} the probability measures \(\cQQ_1\) and \(\QQ_1\) are given by 
 \begin{align}
 \label{eq:ex1}
        &\cQQ_1(E_0 \times E_1) =\PP^{\cpolicy}_{\idist}(\proc_o \in E_0,\proc_1 \in E_1)\nn\\
        &=\int_{E_0}  \idist(\odif{\stt{}{0}})  \PP^{\cpolicy}_{\idist}\cprobof[\big]{\proc_1 \in E_1 \given \stt{}{0}},\nn\\
        &= \int_{E_0}  \idist(\odif{\stt{}{0}}) 
        \PP^{\cpolicy}_{\idist}\Big(\big(a_{11} \stt{1}{0} + b_{1}\ccont{}{0,1} + \pnoise{1}{0}, \cdots \nn\\ &\qquad \cdots, a_{22} \stt{2}{0} + b_{2}\bigl(k_{0,2}\stt{2}{0} + \pe{}{0,2}\bigr) + \pnoise{2}{0}\big) \in E_1\Big),
 \end{align}
 and 
\begin{align}\label{eq:ex11}
        &\QQ_1(E_0 \times E_1) = \PP^{\policy}_{\idist}(\proc_o \in E_0,\proc_1 \in E_1)\nn\\
        &=\int_{E_0}  \idist(\odif{\stt{}{0}})  \PP^{\policy}_{\idist}\cprobof[\big]{\proc_1 \in E_1 \given \stt{}{0}},\nn\\
        &=\int_{E_0}  \idist(\odif{\stt{}{0}}) 
        \PP^{\policy}_{\idist}\Big(\big(a_{11} \stt{1}{0} + b_{1}(k_{0,1} \stt{1}{0} + \pe{}{0,1})+\pnoise{1}{0}, \cdots \nn\\ &\qquad \cdots, a_{22} \stt{2}{0} + b_{2}\bigl(k_{0,2}\stt{2}{0} + \pe{}{0,2}\bigr) + \pnoise{2}{0}\big) \in E_1\Big).
 \end{align}
  Note that in order to show \(\QQ_1 \abscont \cQQ_1\) it is enough to demonstrate that for any set \(E_1 \in \Borelsigalg(\stsp)\) such that \(\PP^{\cpolicy}_{\idist}\cprobof[\big]{\proc_1 \in E_1 \given \stt{}{0}} = 0\) will imply that \(\PP^{\policy}_{\idist}\cprobof[\big]{\proc_1 \in E_1 \given \stt{}{0}} =0.\) Choose \((\pe{}{t})_{t \in \Nz}\) such that \((\pe{}{t,1})_{t \in \Nz}\) is independent of \((\pe{}{t,2})_{t \in \Nz} \) and
 \((\pe{}{t,1})_{t \in \Nz}, (\pe{}{t,2})_{t \in \Nz} \sim \mathsf{Unif}([c_1, c_2])\) where \(c_1, c_2 > 2\abs{c}\) at time \(t\). \textcolor{black}{We observe that if the attack scheme is such that \(\ccont{}{0,1} \sim \mathsf{Unif}([-c_2,-c_1])\), then absolutely continuity of \(\QQ_1\) with respect to \(\cQQ_1\) will not hold,} i.e., for the particular choice of \(\pe{}{t}\), the attacker comes up with a suitable control action \(\ccont{}{t}\) to remain hidden. As a matter of fact, it does not hold for any choice of private excitation. 
\end{example}
The above simple example demonstrates that suitable network structure is a \emph{necessity} for \(\QQ_1 \abscont \cQQ_1\) and hence, \(\QQ_n \abscont \cQQ_n\) for every \(n \in \Nz\) via the principle of mathematical induction. In other words, it is necessary that each agent is influenced either directly or indirectly by at least one of the honest actuators. This condition enables the private excitation to affect all the malicious actuators (if they exists). 

\begin{example}[Choice of the private excitation]\label{exam:support2}
 The objective here is to understand
the role of suitable private excitation in the separability of the state trajectories under \(\policy\) and \(\cpolicy\). In Example \ref{exam:graph2}, let us assume that \[A = \begin{pmatrix}
a_{11} & a_{12}\\
0 & a_{22}
\end{pmatrix}.\]
where \(a_{11},a_{12},a_{22} \neq 0\) and the rest of the data remains identical to Example \ref{exam:graph2}. Firstly, we observe that if the private excitation is chosen as in Example \ref{exam:graph2}, then with the attack scheme \(\ccont{}{0,1} \sim \mathsf{Unif}([-c_2,-c_1])\) adopted by the adversary, \(\QQ_1 \abscont \cQQ_1\) is not guaranteed. On the other hand, if the underlying distribution of the private excitation is fully supported on \(\Rbb\), then \(\QQ_1 \abscont \cQQ_1 \) holds. 
\end{example} 
In short, even though the network structure is such that every malicious agent gets influenced by at least one honest agent, the choice of the private excitations play a pivotal role and must be selected carefully.}

\subsection{Discussion}
\label{subsec:tech_discuss}
Let us discuss several important properties of Definition \ref{def:separator} and our main results --- Theorem \ref{p:key}, Theorem \ref{prop:main}, and Corollary \ref{prop:hmain}.
{\color{black}\begin{rem}
\textbf{(On private excitations)} In this work we employ private excitation to ensure that \(\QQ_n \abscont \cQQ_n\) holds for each \(n \in \Nz\), which suggests that it plays a crucial role in achieving a fully supported probability measure \(\PP^{\cpolicy}_{\idist}\). Moreover, the method developed in this article depends on the nature of the infinitely long sample paths \(\csample\) that allow us to establish sufficient conditions for their non-separability. Our results are therefore, fundamentally different from the ones provided in \cite{ref:YM-BS-09,ref:BS-PRK-16} and \cite{ref:PH-MP-RV-AA-17,ref:PH-MP-RV-AA-18} where watermarking techniques are used to provide a unique signature to the output signal and are only effective when the control actions are hidden.
\end{rem}}
\begin{rem}
\textbf{(On negative and positive results)} The sufficient conditions in Theorem \ref{prop:main} (and  the assertion \ref{it:cor-th2} in Corollary \ref{prop:hmain}) are classified into \emph{positive} result --- \ref{it:ta2}, and \emph{negative} results --- \ref{it:fa2}-\ref{it:faaaa2}, respectively. The \emph{positive result} motivates us to choose a certain class of private excitation with suitable statistics. Indeed, a designer's objective must be to ensure that the Radon-Nikodym derivative in \eqref{e:key} exists for each \(n \in \Nz\) and it converges to zero almost surely. The only apparatus available to them is the private excitation and therefore, one must carefully choose them so that the hypotheses of \ref{it:ta2} hold. On the other hand, the \emph{negative results} can be interpreted as the \emph{fundamental limitations} of the chosen statistics for the private excitation, which is also crucial for a designer. Indeed if the designer has a priori information that the system may be subjected to certain types of malicious attacks so that the generated sample path satisfies the hypotheses of \ref{it:fa2} and \ref{it:faaaa2}, they should refrain from injecting a Gaussian process as their private excitation. The negative results also suggest that it is essential to go \emph{beyond} the regime of Gaussian private excitation to secure the CPS against intelligent and sophisticated attacks, an issue that will be addressed in subsequent articles.
\end{rem}

\begin{rem}
The proofs of Theorem \ref{prop:main} (and Corllary \ref{prop:hmain}) make use of \emph{Jessen's theorem} \cite[Theorem 5.2.20]{ref:Str-11} which gives necessary and sufficient conditions for mutual singularity of two measures \(\PP^{\policy}_{\idist}\) and \(\PP^{\cpolicy}_{\idist}\). Specifically, we check local absolute continuity \cite[Definition 1, p. 165, vol. 2]{ref:Shi-19} of \(\PP^{\policy}_{\idist}\) with respect to \(\PP^{\cpolicy}_{\idist}\), and the convergence of the corresponding Radon-Nikodym derivative (if it exists) \emph{almost surely}. As per our knowledge, results pertaining to almost everywhere convergence are presented only in \cite{ref:BS-PRK-16} and later as an extension in \cite{ref:PH-MP-RV-AA-17}; however, those technical results did not employ Jessen's theorem.
\end{rem}
\begin{rem}
In Theorem \ref{prop:main} (and Corollary \ref{prop:hmain}) the sequences \(\bigl(\stt{}{t}-\meantr{\policy}{t}\bigr)_{t \in \Nz}\) and \(\bigl(\cstt{}{t}-\meantr{\cpolicy}{t}\bigr)_{t \in \Nz}\) denote the \emph{path deviation} of the state trajectories \((\stt{}{t})_{t \in \Nz}\) and \((\cstt{}{t})_{t \in \Nz}\) from their conditional means \(\bigl(\meantr{\policy}{t}\bigr)_{t \in \Nz}\) and \(\bigl(\meantr{\cpolicy}{t}\bigr)_{t \in \Nz}\), respectively. Consequently,
\begin{itemize}
    \item \(\sum_{t=1}^n \eig{\inverse}{\min,t}\norm{\stt{}{t} - \meantr{\policy}{t}}^2\) corresponds to the cumulative weighted energy of the path deviation under the influence of \(\policy\) (i.e., under no-attack condition) until time instant \(t = n\), and
    \item \(\sum_{t=1}^n \ceig{\inverse}{\max,t}\norm{\cstt{}{t} - \meantr{\cpolicy}{t}}^2\) corresponds to the cumulative weighted energy of the path deviation under the influence of \(\cpolicy\) (under attack condition) until time instant \(t = n\).
\end{itemize}
Now let us examine the statements of Theorem \ref{prop:main} in some detail:
\begin{itemize}[label = \(\circ\),leftmargin=*]
    \item The hypothesis of \ref{it:fa2} refers to the fact that if the corrupt policy \(\cpolicy\) is such that the cumulative weighted energy of the path deviations under the influence of \(\policy\) and \(\cpolicy\) gradually become comparable in the limit \(n \ra + \infty\), then it is not possible to separate the corresponding state trajectories \((\stt{}{t})_{t \in \Nz}\) and \((\cstt{}{t})_{t \in \Nz}\). Note that this condition does not require the state trajectories \((\stt{}{t})_{t \in \Nz}\) and \((\cstt{}{t})_{t \in \Nz}\) to be identical across \(t \in \Nz\). 
    
    \item The hypotheses of \ref{it:faaaa2} state that if the ratio of the cumulative weighted energy of the path deviation under the influence of \(\policy\) and \(\cpolicy\), respectively, is uniformly bounded across \(n \in \Nz\), and the cumulative weighted energy of the path deviation under the influence of \(\cpolicy\) is also finite, then we cannot guarantee the existence of a separator in sense of Definition \ref{def:separator}.\footnote{The condition \(\sup_{n \in \Nz}\rate_n \leqslant M\) is reminiscent of the boundedness hypothesis in \emph{Stochastic Approximation} \cite{ref:VSB-08} which has now become standard.}
    
    \item While the assertions of \ref{it:fa2} and \ref{it:faaaa2} are negative results, \ref{it:ta2} is a positive result that asserts the existence of a separator. It states that trajectories \((\stt{}{t})_{t\in \Nz}\) and \((\cstt{}{t})_{t\in \Nz}\) under \(\policy\) and \(\cpolicy\), respectively, can be separated in the sense of Definition \ref{def:separator} if the cumulative weighted energy of the path deviation \((\stt{}{t} - \meantr{\policy}{t})_{t \in \Nz}\) under \(\policy\) increases monotonically over \(n \in \Nz\) at a relatively faster rate than the corresponding cumulative weighted energy of \((\cstt{}{t} - \meantr{\cpolicy}{t})_{t \in \Nz}\) under \(\cpolicy\). 
    A careful scrutiny suggests that the two noteworthy features, namely, asymptotic properties and path properties of the CPS, get reflected in the assertions of Theorem \ref{prop:main}.
\end{itemize}
\end{rem}
\textcolor{black}{\begin{rem}
    \label{rem:on fdi}
    \textbf{(False data injection (FDI) attacks)} A popular type of attack is the false data injection on the actuators, where an intelligently crafted malicious control action \(\ccont{}{t}\) at time \(t\) is typically given by 
    \begin{align*}
        \ccont{}{t} = \cont{}{t} + d_t &= \pmap{t}(\cphist_t) + \pe{}{t} + d_t\\
        &= \cpmap{t}(\cphist_t) + \pe{}{t},
    \end{align*}
	where \(\cpmap{t}(\cphist_t) = \pmap{t}(\cphist_t) + d_t\), \(\bigl(\cphist_t\bigr)_{t \in \Nz}\) are generated under \(\cpolicy\), and \(d_t\) denotes the adversarial disturbance injected by the attacker into the communication channel (see Fig. \ref{fig:attack model}). Here \(\pe{}{0} \sim \grv(0,\var_{\pe{}{}})\) and Assumption \ref{assum:gaussian case} is enforced. We observe that such attacks are history dependent and the framework developed in this article naturally subsumes these attacks. Moreover, a careful scrutiny reveals that
    \begin{itemize}[label = \(\circ\),leftmargin=*]
        \item  under \(\cpolicy_t\) at time \(t\), \(\proc_t\) is a conditionally Gaussian random vector given \(\cphist_t\), with the conditional mean
        \begin{align*}
            \meantr{\cpolicy}{t+1}  &= \EE^{\cpolicy}_{\idist} \cexpecof[\big]{\proc_{t+1} \given \sigalg_{t}}(\csample) = \EE^{\cpolicy}_{\idist} \cexpecof[\big]{\proc_{t+1} \given \cphist_t}
            \nn\\&= A\cstt{}{t}+B \EE^{\cpolicy}_{\idist} \cexpecof[\big]{\control_t \given \cphist_t}+\EE^{\cpolicy}_{\idist} \cexpecof[\big]{\pnrv{t} \given \cphist_t}, \\
            &= A\cstt{}{t}+B \begin{pmatrix}
                \cpmap{t,1}(\cphist_t)\\
                \pmap{t,2}(\cphist_t)   
            \end{pmatrix},
        \end{align*}
        and the conditional variance
        \begin{align*}
            \var_{t+1}^{\cpolicy} &= B \tilde{\var}B^{\top} + \var_{\pnrv{}},
        \end{align*}
        where \(\tilde{\var} = \EE\cexpecof[\big]{\pe{}{t}\pe{\top}{t} \given \cphist_t} = \var_{\pe{}{}}\);
        \item under \(\policy_t\) at time \(t\), \(\proc_t\) is a conditionally Gaussian random vector given \(\phist_t\), with the conditional mean  
        \(
            \meantr{\policy}{t+1} = A\stt{}{t} + B \pmap{t}(\phist_t)   
        \)
        and the conditional variance being 
        \(
            \var^{\policy}_{t}  = B \var_{\pe{}{}}B^{\top} + \var_{\pnrv{}}.    
        \)
    \end{itemize}
The aforementioned two observations indicate that Lemma \ref{lem:rvg} holds and the rest of the proofs for Theorem \ref{prop:main} (and Corollary \ref{prop:hmain}) remain unchanged in this special case of FDI attacks.
\end{rem}
}

\section{Conclusions and future work}
\label{sec:conclusion}
{\color{black} This article established a general framework for almost sure detectability of arbitrarily clever actuator attacks. To that end, a novel idea of the separability of state trajectories was introduced. It was demonstrated that the idea of separability can be leveraged in the context of security of CPS for detecting the presence of malicious actuators in them; Theorem \ref{p:key} is the main foundational result. For a stochastic linear time-invariant system, the mechanism of injecting private excitations was adopted for the purpose of almost sure detection of the presence of malicious attacks, and sufficient conditions for the existence and non-existence of separators were established.

Although the assertions established in this article are based on strong guarantees (almost sure detectability) and remain valid for a vast class of intelligent actuator attacks, these results are asymptotic, and in general difficult to verify. We emphasize that almost sure detection is not possible in finite time due to the canonical structure of the underlying probability space. However, the fundamental observations developed in this article naturally produce applicable detection (with high probabilistic guarantees) algorithms along the lines of statistical hypothesis testing; these extensions will be reported in subsequent articles.}

 \appendix
\section*{Appendix}

\section{Existence of the measure \(\PP^{\policy}_{\idist}\) and \(\PP^{\cpolicy}_{\idist}\)}
\label{appen:eom}
This appendix includes a proof of the existence of \(\PP^{\policy}_{\idist}\) and \(\PP^{\cpolicy}_{\idist}\) corresponding to \(\policy\) and \(\cpolicy\), respectively. {\color{black}Below we outline the steps involved:
\begin{itemize}[label = \(\circ\)]
    \item Restrict \(\tRR^{\policy}_{\idist}\) and \(\tRR^{\cpolicy}_{\idist}\) to the finite-dimensional measurable space \((\ps,\tsigalg_n)\) where \(n \in \Nz\), namely \(\tQQ_n\) and \(\ctQQ_n\);
    \item Marginalize \(\tQQ_n\) and \(\ctQQ_n\) with respect to \((\cont{}{t})_{t \in \Nz}\) and \((\ccont{}{t})_{t \in \Nz}\), respectively;
    \item Extend \(\tQQ_n\) and \(\ctQQ_n\) to \(\bigl(\ps,\Borelsigalg(\ps)\bigr)\) by Kolmogorov's extension theorem \cite[Theorem 3, p. 197]{ref:Shi-16}.
\end{itemize}
}

Recall the filtered probability space  \(\Bigl(\tss,\Borelsigalg{(\tss),\bigl(\tsigalg_t \bigr)_{t \in \Nz}},\tRR^{\policy}_{\idist} \Bigr).\) Let us define the measure \(\tQQ_n\) on the measurable space \((\tss,\tsigalg_n)\) by
\begin{align*}
    &\tQQ_n(E) \\&= \int_{S_0} \idist(\odif{y_0}) \int_{T_0} \policy_0(\hist_0; \odif{\cont{}{0}})\int_{S_1} \tRR^{\policy}_{\idist} \cprobof[\big]{\proc_1 \in \odif{y_1} \given \hist_0, \cont{}{0}}  \cdots \\
	& \qquad \cdots  \int_{T_{n-1}} \policy_{n-1}(\hist_{n-1}; \odif{\cont{}{n-1}}) \int_{S_{n}} \tRR^{\policy}_{\idist} \cprobof[\big]{\proc_n \in \odif{y_n} \given \hist_{n-1}, \cont{}{n-1}} \nn
\end{align*}
where \(E \in \tsigalg_n\) is of the form \(E = S_0 \times T_0 \times S_1 \times T_1 \times \cdots \times S_n,\) with the Borel sets \(S_0,S_1,\ldots, S_n \in \Borelsigalg(\stsp)\) and \(T_0,T_1,\ldots, T_{n-1} \in \Borelsigalg(\admact).\) Then, for cylindrical sets of the form \(E' = S_0 \times S_1 \times \cdots \times S_n\) where \(S_i \in \Borelsigalg{(\stsp)}\) for \(i=1,\ldots,n\) we define the measure of \(E'\) by 
\begin{align}
\label{appen-eq:emarginal}
    &\QQ_n (E') \nn\\&\Let \int_{E'}\biggl[\int_{\underbrace{\admact \times \cdots \times \admact}_{\text{n-times}}} \tQQ_n(\odif{\stt{}{0}},\odif{\cont{}{0}},\odif{\stt{}{1}},\ldots,\odif{\cont{}{n-1}},\odif{\stt{}{n}})\biggr] \nn\\
    &\stackrel{\mathclap{(\star)}}{=} \int_{E'}\biggl[\underbrace{\int_{\admact}  \cdots \int_{\admact}}_{\text{n-times}} \idist(\odif{y_0})  \policy_0(\hist_0; \odif{\cont{}{0}}) \tRR^{\policy}_{\idist} \cprobof[\big]{\proc_1 \in \odif{y_1} \given \hist_0, \cont{}{0}}  \cdots \nn\\
	& \hspace{8mm} \qquad \cdots   \policy_{n-1}(\hist_{n-1}; \odif{\cont{}{n-1}})  \tRR^{\policy}_{\idist} \cprobof[\big]{\proc_n \in \odif{y_n} \given \hist_{n-1}, \cont{}{n-1}}                \biggr] \nn \\
	&= \int_{S_0}\idist(\odif{y_0}) \int_{\admact} \policy_0(\hist_0; \odif{\cont{}{0}}) \int_{S_1}\tRR^{\policy}_{\idist} \cprobof[\big]{\proc_1 \in \odif{y_1} \given \hist_0, \cont{}{0}}  \cdots \nn\\
	& \quad  \cdots   \int_{\admact}\policy_{n-1}(\hist_{n-1}; \odif{\cont{}{n-1}}) \int_{S_n} \tRR^{\policy}_{\idist} \cprobof[\big]{\proc_n \in \odif{y_n} \given \hist_{n-1}, \cont{}{n-1}}, \tag{(A.1)}
\end{align}
where \((\star)\) follows from Fubini's theorem \cite[Theorem 8, p. 235]{ref:Shi-16}. Observe that for each \(n \in \Nz\), \(\QQ_n\) is a measure on \((\stsp^n, \Borelsigalg(\stsp^n))\). Clearly, \(\QQ_{n+1}(E' \times \stsp) = \QQ_n(E')\); invoke Kolmogorov's extension theorem \cite[Theorem 3, p. 197]{ref:Shi-16} to assert the existence of a unique measure \(\PP^{\policy}_{\idist}\) on \((\stsp^{\Nz}, \Borelsigalg{(\stsp^{\Nz})})\) such that  
\(
\PP^{\policy}_{\idist} \probof[\big]{\proc_0 \in S_0,\proc_1 \in S_1, \ldots, \proc_{n} \in S_n} = \QQ_n(E').
\)
The proof for the existence of the measure \(\PP^{\cpolicy}_{\idist}\) induced by \(\cpolicy\) proceeds along similar lines.
\hfill \(\blacksquare\)
\section{Proofs}
\label{appen:proof_auxx}
Here we establish the proof of Lemma \ref{lem:wrpm}
\begin{pf}\label{appen:pa}
From \eqref{appen-eq:emarginal}, we start by observing that for every \(n \in \Nz\),
\begin{align*}
        &\int_{\admact}\policy_{n-1}\cprobof[\big]{\hist_{n-1};\odif{\cont{}{n-1}}} \int_{S_n}\tRR^{\policy}_{\idist} \cprobof[\big]{\proc_{n} \in \odif{y_{n}}\given y_{n-1},\cont{}{n-1}}\\
        &\stackrel{\mathclap{(1)}}{=} \int_{\admact}\policy_{n-1}\cprobof[\big]{\hist_{n-1};\odif{\cont{}{n-1}}} \tRR^{\policy}_{\idist} \cprobof[\big]{\proc_{n} \in S_n\given y_{n-1},\cont{}{n-1}}\\
        &\stackrel{\mathclap{(2)}}{=} \EE^{\policy}_{\idist} \cexpecof[\Big]{\tRR^{\policy}_{\idist} \cprobof[\big]{\proc_{n} \in S_n\given \hist_{n-1},\cont{}{n-1}} \given \hist_{n-1}}\\
        & = \tRR^{\policy}_{\idist} \cprobof[\big]{\proc_{n} \in S_n\given \hist_{n-1}}
\end{align*}
where \((1)\) and \((2)\) above follows from \eqref{eq:updf}. From \eqref{e:marginal},
\begin{align}\label{eq:auxeq}
	&\QQ_n(E) \nn\\& \Let \int_{S_0} \idist(\odif{y_0})  \int_{\admact} \policy_0(\hist_0; \odif{\cont{}{0}})\int_{S_1} \tRR^{\policy}_{\idist} \cprobof[\big]{\proc_{1} \in \odif{y_{1}}\given y_{0},\cont{}{0}}  \cdots \nn\\
	& \quad \cdots  \int_{\admact} \policy_{n-1}(\hist_{n-1}; \odif{\cont{}{n-1}}) \int_{S_{n}} \tRR^{\policy}_{\idist} \cprobof[\big]{\proc_{n} \in \odif{y_{n}}\given y_{n-1},\cont{}{n-1}} \nn\\
	&= \int_{S_0} \idist(\odif{y_0})  \int_{\admact} \policy_0(\hist_0; \odif{\cont{}{0}})\int_{S_1} \tRR^{\policy}_{\idist} \cprobof[\big]{\proc_{1} \in \odif{y_{1}}\given y_{0},\cont{}{0}}  \cdots \nn\\
	& \qquad \cdots \int_{S_n} \tRR^{\policy}_{\idist} \cprobof[\big]{\proc_{n} \in \odif{y_n}\given \hist_{n-1}} \tag{(B.1)}
\end{align}
Also, note that 
\begin{align*}
    &\int_{\admact}\policy_{n-2}\cprobof[\big]{\hist_{n-2};\odif{\cont{}{n-2}}} \int_{S_{n-1}}\tRR^{\policy}_{\idist} \cprobof[\big]{\proc_{n-1} \in \odif{y_{n-1}}\given y_{n-2},\cont{}{n-2}} \cdots\\ & \qquad \qquad \qquad \qquad \qquad \qquad \cdots \int_{S_n} \tRR^{\policy}_{\idist} \cprobof[\big]{\proc_{n} \in \odif{y_n}\given \hist_{n-1}}\\
    &= \int_{\admact}\policy_{n-2}\cprobof[\big]{\hist_{n-2};\odif{\cont{}{n-2}}} \tRR^{\policy}_{\idist} \cprobof[\big]{\proc_{n-1} \in S_{n-1}, \proc_{n} \in S_{n}\given \hist_{n-2}, \cont{}{n-2}}\\
    &= \tRR^{\policy}_{\idist} \cprobof[\big]{\proc_{n-1} \in S_{n-1}, \proc_{n} \in S_{n}\given \hist_{n-2}}\\
    &= \int_{S_{n-1}}\tRR^{\policy}_{\idist} \cprobof[\big]{\proc_{n-1} \in \odif{y_{n-1}}\given \hist_{n-2}}\int_{S_n} \tRR^{\policy}_{\idist} \cprobof[\big]{\proc_{n} \in \odif{y_n}\given \hist_{n-2},y_{n-1}},
\end{align*}
which implies that 
\begin{align*}
&\QQ_n(E)  \\&\Let \int_{S_0} \idist(\odif{y_0})  \int_{\admact} \policy_0(\hist_0; \odif{a_0})   \int_{S_1} \tRR^{\policy}_{\idist} \cprobof[\big]{\proc_{1} \in \odif{y_{1}}\given y_{0},\cont{}{0}} \cdot\cdot \\ & \cdot\cdot\int_{S_{n-1}}\tRR^{\policy}_{\idist} \cprobof[\big]{\proc_{n-1} \in \odif{y_{n-1}}\given \hist_{n-2}}  \int_{S_n} \tRR^{\policy}_{\idist} \cprobof[\big]{\proc_{n} \in \odif{y_n}\given \hist_{n-2},y_{n-1}}.
\end{align*}
Using \eqref{eq:auxeq} and proceeding backward in a similar way, we obtain
\begin{align*}
\QQ_n(E)  \Let \int_{S_0} \idist(\odif{y_0})  &\int_{S_1} \tRR^{\policy}_{\idist} \cprobof[\big]{\proc_{1} \in \odif{y_{1}}\given \phist_0}  \cdots \\ &\cdots\int_{S_n} \tRR^{\policy}_{\idist} \cprobof[\big]{\proc_{n} \in \odif{y_n}\given \phist_{n-1}}, 
\end{align*}
where \(\phist_t = (\stt{}{0},\ldots, \stt{}{t})\) and the first assertion follows immediately by observing that \(\tRR^{\policy}_{\idist}\probof[\big]{\proc_{0} \in S_0, \cont{}{0} \in \admact, \proc_{1} \in S_1, \cont{}{1} \in \admact, \ldots,\proc_n \in S_n} = \PP^{\policy}_{\idist}\probof[\big]{\proc_{0} \in S_0, \proc_{1} \in S_1, \ldots,\proc_n \in S_n}\). Also, note that 
\begin{align*}
    &\int_E \QQ_n\probof[\big]{\odif{y_0},\ldots, \odif{y_n}}\\
    &= \int_{S_0} \idist(\odif{y_o}) \prod_{t=1}^n \int_{S_{t}} \PP^{\policy}_{\idist} \cprobof[\big]{\proc_t \in \odif{y_t} \given \phist_{t-1}},
\end{align*}
which proves the second assertion. \hfill \(\blacksquare\)
\end{pf}
\subsection{Proof of Theorem \ref{prop:main}} \label{appen:main_proof}
Here we establish some auxiliary results for the proof of Theorem \ref{prop:main}. We start by giving a proof for Lemma \ref{lem:rvg}.
{\color{black}\begin{pf}[Proof of Lemma \ref{lem:rvg}]
Recall that by hypothesis, all policies are in \(\pclass{RH}\). Firstly, we consider the case where the CPS \eqref{eq:lsyscompactt} has at least one malicious actuator. Given the history \(\cphist_t\), under the corrupted policy \((\cpolicy_t)_{t \in \Nz}\), we have \(\proc_{t+1}=A\cstt{}{t}+B\control_t + W_t\) where \(W_t\) denotes the process noise at \(t \in \Nz\) and \(U \Let (\control_t)_{t \in \Nz}\) denotes the control process. Let \(\meantr{\cpolicy}{t}\) and \(\var_t^{\cpolicy}\) denotes the conditional mean and variance of \(\proc_t\) given the history \(\cphist_{t-1}\). We observe that \(\proc_{t+1}\) is the sum of two Gaussian random vectors \(\control_{t}\) and \(W_t\). Hence, \(\proc_t\) is a conditional Gaussian random vector with the conditional mean being
\begin{align}\label{eq:c_cmean}
    \meantr{\cpolicy}{t+1}  &= \EE^{\cpolicy}_{\idist} \cexpecof[\big]{\proc_{t+1} \given \sigalg_{t}}(\csample) = \EE^{\cpolicy}_{\idist} \cexpecof[\big]{\proc_{t+1} \given \cphist_t}
    \nn\\&= A\cstt{}{t}+B \EE^{\cpolicy}_{\idist} \cexpecof[\big]{\control_t \given \cphist_t}+\EE^{\cpolicy}_{\idist} \cexpecof[\big]{\pnrv{t} \given \cphist_t}, \tag{(B.3)}
\end{align} 
and the conditional variance being 
\begin{equation}\label{eq:c_cvar}
    \var_{t+1}^{\cpolicy} = B \tilde{\var}B^{\top} + \var_{\pnrv{}}. \tag{(B.4)}
\end{equation} Here 
\begin{align*}
\EE^{\cpolicy}_{\idist} \cexpecof[\big]{\control_{t} \given \cphist_t} &= \bigl(\EE^{\cpolicy}_{\idist} \cexpecof[\big]{\cpmap{t,1}(\cphist_t) \given \cphist_t}, \EE^{\cpolicy}_{\idist} \cexpecof[\big]{\pmap{t,2}(\cphist_t) + \pe{}{t,2}\given \cphist_t}\bigr) \\&= \bigl(\cpmap{t,1}(\cphist_t),\pmap{t,2}(\cphist_t) \bigr),
\end{align*}
and 
 \begin{align*} 
    \tilde{\var} &= \EE^{\cpolicy}_{\idist} \cexpecof[\Big]{\Bigl(\control_{t} - \EE^{\cpolicy}_{\idist} \cexpecof[\big]{\control_{t}\given \cphist_t} \Bigr)\Bigl(\control_{t} - \EE^{\cpolicy}_{\idist} \cexpecof[\big]{\control_{t}\given \cphist_t} \Bigr)^{\top}\given \cphist_t}\\
    &= \begin{pmatrix}
    \zv{M} && \zv{M \times \node - M}\\
    \zv{\node - M \times M} && \var_2
    \end{pmatrix}.
 \end{align*}
Next, we consider the case when all the actuators are honest. Arguments similar to the preceding case hold when an honest policy
\((\policy_t)_t\) is employed:  given the history \(\phist_t\), \(\proc_t\) is the sum of two Gaussian random vectors \(\control_{t}\) and \(W_t\), consequently, \(\proc_{t+1}\) is also a conditional Gaussian random vector with the  conditional mean being
\begin{align}\label{eq:c_mean}
\meantr{\policy}{t+1} &=\EE^{\policy}_{\idist} \cexpecof[\big]{\proc_{t+1} \given \sigalg_{t}}(\sample) = \EE^{\policy}_{\idist}\cexpecof[\big]{\proc_{t+1} \given \phist_t} \nn\\
&= A\stt{}{t} + B\EE^{\policy}_{\idist} \cexpecof[\big]{\control_{t}\given \phist_t} + \EE^{\policy}_{\idist} \cexpecof[\big]{\pnrv{t} \given \phist_t} \nn\\
 &= A\stt{}{t} + B\EE^{\policy}_{\idist} \cexpecof[\big]{\pmap{t}(\phist_t) + \pe{}{t}\given \phist_t} + \EE^{\policy}_{\idist} \cexpecof[\big]{\pnrv{t} \given \phist_t}\nn\\
&= A\stt{}{t} + B \pmap{t}(\phist_t), \tag{(B.5)}
\end{align}
and the conditional variance being
\begin{equation}\label{eq:c_var}
    \var^{\policy}_{t}  = B \var_{\pe{}{}}B^{\top} + \var_{\pnrv{}}. \tag{(B.6)}
\end{equation} This completes the proof. \hfill \(\blacksquare\)
\end{pf}
 The next lemma establishes absolute continuity of the measure \(\QQ_n \) with respect to \(\cQQ_n\) for each \(n\in \Nz\).\footnote{See ``Notation'' in \S\ref{sec:intro} for the definition of two measures being absolutely continuous with respect to each other, and Appendix \ref{appen:prob basics} for a brief description on some of the terms used in the subsequent results.} 
\begin{lem}
   \label{lem:eRND}
   If the hypotheses of Theorem \ref{prop:main} hold, then \(\QQ_n \abscont \cQQ_n\) for each \(n \in \Nz\).
\end{lem}
\begin{pf}
The proof proceeds by the principle of mathematical induction (we shall adhere to various notations in Theorem \ref{prop:main}). Our base case is \(n=1\), and from Lemma \ref{lem:rvg} it follows that \(\proc_1\) is a conditional Gaussian random vector given \(\proc_0\) under the influence of the policies \((\policy_t)_{t \in \Nz}\) and \((\cpolicy_t)_{t \in \Nz}\). Observe that from \eqref{eq:lem:wrpm},
\begin{align*}
        &\QQ_1 \probof[\big]{\odif{\stt{}{0}}, \odif{\stt{}{1}}} \\&= \idist(\odif{\stt{}{0}}) \PP^{\policy}_{\idist} \cprobof[\big]{\proc_1 \in \odif{\stt{}{1}} \given \phist_0}\\
        &\stackrel{\mathclap{(\ast)}}{=} \idist(\odif{\stt{}{0}}) \frac{\PP^{\policy}_{\idist} \cprobof[\big]{\proc_1 \in \odif{\stt{}{1}} \given \phist_0}}{\PP^{\cpolicy}_{\idist} \cprobof[\big]{\proc_1 \in \odif{\stt{}{1}} \given \phist_0}} \PP^{\cpolicy}_{\idist} \cprobof[\big]{\proc_1 \in \odif{\stt{}{1}} \given \phist_0},
\end{align*}
where the term on the right-hand side of \((\ast)\) is well-defined and can be expanded as 
\begin{align*}
        &\frac{\PP^{\policy}_{\idist} \cprobof[\big]{\proc_1 \in \odif{\stt{}{1}} \given \phist_0}}{\PP^{\cpolicy}_{\idist} \cprobof[\big]{\proc_1 \in \odif{\stt{}{1}} \given \phist_0}} \\&= \frac{\det\bigl(\var^{\cpolicy}_{1}\bigr)^{1/2}}{\det\bigl(\var^{\policy}_{1}\bigr)^{1/2}} \exp\Biggl[-\frac{1}{2} \Bigl(\stt{}{1} - \meantr{\policy}{1} \Bigr)^{\top} \Bigl(\var^{\policy}_{1} \Bigr)^{-1}\Bigl(\stt{}{1} - \meantr{\policy}{1} \Bigr) \\ & \qquad + \frac{1}{2} \Bigl(\stt{}{1} - \meantr{\cpolicy}{1} \Bigr)^{\top} \Bigl(\var^{\cpolicy}_{1} \Bigr)^{-1}\Bigl(\stt{}{1} - \meantr{\cpolicy}{1} \Bigr) \Biggr],
\end{align*}
where \(\meantr{\policy}{1}\) and \(\meantr{\policy}{1}\) are defined in \eqref{eq:c_mean} and \eqref{eq:c_cmean}, respectively. Here the quantities \(\var^{\policy}_{1}\) and \(\var^{\cpolicy}_{1}\) denote the conditional variances of \(\proc_1\) given \(\proc_0\), under the policy \((\policy_t)_{t \in \Nz}\) and \((\cpolicy_t)_{t \in \Nz}\) respectively. This implies that for any \(E \in \sigalg_1\) of the form \(E = S_0 \times  S_1\) where \(S_0,\,S_1 \in \Borelsigalg(\stsp)\), we have
\begin{align*}
        \QQ_1 \probof[\big]{E} &= \int_E \QQ_1 \probof[\big]{\odif{\stt{}{0}}, \odif{\stt{}{1}}}\\
        &= \int_{E} \idist(\odif{\stt{}{0}}) \frac{\PP^{\policy}_{\idist} \cprobof[\big]{\proc_1 \in \odif{\stt{}{1}} \given \phist_0}}{\PP^{\cpolicy}_{\idist} \cprobof[\big]{\proc_1 \in \odif{\stt{}{1}} \given \phist_0}} \PP^{\cpolicy}_{\idist} \cprobof[\big]{\proc_1 \in \odif{\stt{}{1}} \given \phist_0}\\
        &= \int_E \frac{\PP^{\policy}_{\idist} \cprobof[\big]{\proc_1 \in \odif{\stt{}{1}} \given \phist_0}}{\PP^{\cpolicy}_{\idist} \cprobof[\big]{\proc_1 \in \odif{\stt{}{1}} \given \phist_0}} \idist(\odif{\stt{}{0}})  \PP^{\cpolicy}_{\idist} \cprobof[\big]{\proc_1 \in \odif{\stt{}{1}} \given \phist_0}\\
        &= \int_E \frac{\PP^{\policy}_{\idist} \cprobof[\big]{\proc_1 \in \odif{\stt{}{1}} \given \phist_0}}{\PP^{\cpolicy}_{\idist} \cprobof[\big]{\proc_1 \in \odif{\stt{}{1}} \given \phist_0}} \cQQ_1 \probof[\big]{\odif{\stt{}{0}},\odif{\stt{}{1}}};
\end{align*}
that is, for any \(E \in \sigalg_1\) of the above form, \(\cQQ_1 \probof[\big]{E} = 0\) implies \(\QQ_1 \probof[\big]{E}=0\).  
This proves that the assertion holds for \(n=1\). Suppose now that the assertion holds for an arbitrary but fixed \(n \in \Nz\). Then we have to show that \(\QQ_{n+1} \abscont \cQQ_{n+1}\) in our induction step. To that end, from \eqref{eq:lem:wrpm} we write
\begin{align*}
        &\QQ_{n+1} \probof[\big]{\odif{\stt{}{0}},\ldots,\odif{\stt{}{n}},\odif{\stt{}{n+1}}}\\&=\idist(\odif{y_o}) \prod_{t=1}^{n} \PP^{\policy}_{\idist} \cprobof[\big]{\proc_t \in \odif{y_t} \given \phist_{t-1}} \PP^{\policy}_{\idist} \cprobof[\big]{\proc_{n+1} \in \odif{y_{n+1}} \given \phist_{n}}\\ &=\QQ_n \probof[\big]{\odif{\stt{}{0}},\ldots,\odif{\stt{}{n}}} \PP^{\policy}_{\idist} \cprobof[\big]{\proc_{n+1} \in \odif{\stt{}{n+1}} \given \phist_n}\\
        &\stackrel{\mathclap{(\dagger)}}{=} \prod_{t=1}^n \frac{\PP^{\policy}_{\idist}\cprobof[\big]{\proc_t \in \odif{\stt{}{t}} \given \phist_{t-1}}}{\PP^{\cpolicy}_{\idist}\cprobof[\big]{\proc_t \in \odif{\stt{}{t}} \given \phist_{t-1}}} \cQQ_n \probof[\big]{\odif{\stt{}{0}},\ldots,\odif{\stt{}{n}}} \cdots \\ & \qquad\qquad \qquad \qquad \cdots\PP^{\policy}_{\idist} \cprobof[\big]{\proc_{n+1} \in \odif{\stt{}{n+1}} \given \phist_{n}}\\
        &\stackrel{\mathclap{(\dagger\dagger)}}{=} \,\prod_{t=1}^{n+1} \frac{\PP^{\policy}_{\idist}\cprobof[\big]{\proc_t \in \odif{\stt{}{t}} \given \phist_{t-1}}}{\PP^{\cpolicy}_{\idist}\cprobof[\big]{\proc_t \in \odif{\stt{}{t}} \given \phist_{t-1}}}\cQQ_n \probof[\big]{\odif{\stt{}{0}},\ldots,\odif{\stt{}{n}},\odif{\stt{}{n+1}}}, 
\end{align*}
where the term on the right-hand side of \((\dagger\dagger)\) is well-defined and admits a closed form expression. The equality \((\dagger)\) follows from the induction step for an arbitrary but fixed \(n \in \Nz\) and \((\dagger\dagger)\) follows from \eqref{eq:lem:wrpm}. Choosing \(E \in \sigalg_{n+1}\) of the form \(E = S_0 \times S_1 \times \cdots \times S_{n+1}\) where \(S_0,\ldots, S_{n+1} \in \Borelsigalg(\stsp)\). We integrate 
\begin{align*}
    &\QQ_{n+1}\probof[\big]{\odif{\stt{}{0}},\ldots,\odif{\stt{}{n+1}}}\\&=\prod_{t=1}^{n+1} \frac{\PP^{\policy}_{\idist}\cprobof[\big]{\proc_t \in \odif{\stt{}{t}} \given \phist_{t-1}}}{\PP^{\cpolicy}_{\idist}\cprobof[\big]{\proc_t \in \odif{\stt{}{t}} \given \phist_{t-1}}}\cQQ_n \probof[\big]{\odif{\stt{}{0}},\ldots,\odif{\stt{}{n+1}}},
\end{align*}
on both sides over the set \(E\) to arrive at 
\[
\QQ_{n+1}(E) = \int_E  \prod_{t=1}^{n+1} \frac{\PP^{\policy}_{\idist}\cprobof[\big]{\proc_t \in \odif{\stt{}{t}} \given \phist_{t-1}}}{\PP^{\cpolicy}_{\idist}\cprobof[\big]{\proc_t \in \odif{\stt{}{t}} \given \phist_{t-1}}} \cQQ_{n+1}(\odif{\stt{}{0}},\ldots,\odif{\stt{}{n+1}}).
\]
Arguments similar to the base case \(n=1\) apply here, which proves the assertion for the \((n+1)^{\text{th}}\) induction step. This concludes the proof. 
\hfill \(\blacksquare\)
\end{pf}
We now proceed to the proof of Theorem \ref{prop:main}.
\begin{pf}[Proof of Theorem \ref{prop:main}] 
\textcolor{black}{Recall that we stack the control actions generated by the malicious agents and the control actions generated by the honest agents separately:} for each \(t \in \Nz\), the control actions corresponding to the corrupt policy are partitioned as 
\(
\ccont{}{t} = \bigl(
\ccont{}{t,1}, \cont{}{t,2}
\bigr)
\)
where \(\ccont{}{t,1} = \Bigl(\ccont{1}{t},\cdots,\ccont{M}{t}\Bigr)\) and \(\cont{}{t,2} = \Bigl(\cont{M+1}{t},\cdots,\cont{N}{t}\Bigr)\). Consequently, the input matrix \(B\) and the covariance matrices \(\var_{\pe{}{}}\) and \(\var_{\pnoise{}{}}\) are also partitioned as \(B = \mathsf{diag}\bigl(B_{M},B_{H} \bigr)\) where \(B_{M} \in \Rbb^{M \times M}\) and \(B_H \in \Rbb^{(\node - M) \times (\node -M)}\). Moreover, the variance of the private excitation and the process noise is given by
\begin{align*}
    \var_{\pe{}{}} = \begin{pmatrix}
    \var_1 && \zv{M \times N-M}\\
    \zv{\node - M \times M} && \var_2
    \end{pmatrix}, \quad \var_{\pnoise{}{}} &= \begin{pmatrix}
    \wvar_1 && \wvar_{12}\\
    \wvar_{21} && \wvar_2
    \end{pmatrix}.
\end{align*}
Recall from \eqref{eq:c_var} and \eqref{eq:c_cvar} that the variance of \(\stt{}{t}\) at time \(t\) is given by \[
    \var^{\policy}_{t} = \begin{pmatrix}
        B_M \var_1 B_M^{\top} + \wvar_1 & \wvar_{12}\\
        \wvar_{21} & B_H \var_2 B_H^{\top} + \wvar_2
    \end{pmatrix}
\] 
and the variance of \(\cstt{}{t}\) at time \(t\) is given by \(\var^{\cpolicy}_{t}\)
\[
    \var^{\cpolicy}_{t} = \begin{pmatrix}
         \wvar_1 & \wvar_{12}\\
        \wvar_{21} & B_H \var_2 B_H^{\top} + \wvar_2
    \end{pmatrix}.
\] 
The proof of Theorem \ref{prop:main} is based on the fact that under Assumption \ref{assum:gaussian case}, the random vector \(\proc_t\) is conditionally Gaussian given its entire past \(\phist_{t-1}\) (see Lemma \ref{lem:rvg} for the proof). Moreover, Lemma \ref{lem:eRND} asserts that \(\QQ_n \abscont \cQQ_n\) across all \(n \in \Nz\). This implies that for every \(n \in \Nz\) the Radon-Nikodym derivative \(\frac{\odif{\QQ_n}}{\odif{\cQQ_n}}\) exists \(\PP^{\cpolicy}_{\idist}\mbox{-}\)\emph{almost surely} and is given by the expression
\begin{align}
    \label{eq:RND e}
    \rnd{n}{(\stt{}{0},\ldots,\stt{}{n})} &\Let  \frac{\odif{\QQ_n}}{\odif{\cQQ_n}}(\stt{}{0},\ldots,\stt{}{n})\nn\\
    &= \prod_{t=1}^n \frac{\PP^{\policy}_{\idist}\cprobof[\big]{\proc_t \in \odif{\stt{}{t}} \given \phist_{t-1}}}{\PP^{\cpolicy}_{\idist}\cprobof[\big]{\proc_t \in \odif{\stt{}{t}} \given \phist_{t-1}}}, \tag{(B.7)}
\end{align}
where the term
\begin{align*}
        &\frac{\PP^{\policy}_{\idist}\cprobof[\big]{\proc_t \in \odif{\stt{}{t}} \given \phist_{t-1}}}{\PP^{\cpolicy}_{\idist}\cprobof[\big]{\proc_t \in \odif{\stt{}{t}} \given \phist_{t-1}}}\\ &=  \frac{\det\bigl(\var^{\cpolicy}_{t}\bigr)^{1/2}}{\det\bigl(\var^{\policy}_{t}\bigr)^{1/2}} \,\text{exp}\,\Biggl[-\frac{1}{2} \Bigl(\stt{}{t} - \meantr{\policy}{t} \Bigr)^{\top} \Bigl(\var^{\policy}_{t} \Bigr)^{-1}\Bigl(\stt{}{t} - \meantr{\policy}{t} \Bigr) \\ & \qquad \qquad \qquad \qquad + \frac{1}{2} \Bigl(\stt{}{t} - \meantr{\cpolicy}{t} \Bigr)^{\top} \Bigl(\var^{\cpolicy}_{t} \Bigr)^{-1}\Bigl(\stt{}{t} - \meantr{\cpolicy}{t} \Bigr) \Biggr].
    \end{align*}
This implies that 
\begin{align*}
    &\rnd{n}{(\stt{}{0},\ldots,\stt{}{n})}\\ &\stackrel{\mathclap{(\star)}}{\leqslant} \prod_{t=1}^n \frac{\det\bigl(\var^{\cpolicy}_{t}\bigr)^{1/2}}{\det\bigl(\var^{\policy}_{t}\bigr)^{1/2}} \, \mathrm{exp} \,\Bigg[-\frac{1}{2} \eig{-1}{\min,t}\norm{\stt{}{t} - \meantr{\policy}{t}}^{2} \\& \qquad \qquad \qquad \qquad\qquad+ \frac{1}{2} \ceig{-1}{\max,t} \norm{\stt{}{t} - \meantr{\cpolicy}{t}}^2 \Bigg] \\
    &= \prod_{t=1}^n \frac{\det\bigl(\var^{\cpolicy}_{t}\bigr)^{1/2}}{\det\bigl(\var^{\policy}_{t}\bigr)^{1/2}}  \exp{\Biggl[-\frac{1}{2} \sum_{t=1}^n s_t  +  \frac{1}{2} \sum_{t=1}^n  \breve{s_t}\Biggr]}\\
    &= \prod_{t=1}^n \frac{\det\bigl(\var^{\cpolicy}_{t}\bigr)^{1/2}}{\det\bigl(\var^{\policy}_{t}\bigr)^{1/2}}  \exp{\Biggl[-\frac{1}{2}\Bigl(\rate_n - 1\Bigr) \sum_{t=1}^n  \breve{s_t}\Biggr]}, 
\end{align*}
where \(\rate_n\) is defined in \eqref{eq:rate} and we abbreviate \(s_t \Let\eig{-1}{\min,t} \norm{\stt{}{t} - \meantr{\policy}{t}}^2\) and \(\breve{s_t} \Let \ceig{-1}{\max,t} \norm{\stt{}{t} - \meantr{\cpolicy}{t}}^2\); \(\eig{}{\min,t}\) and \(\ceig{}{\max,t}\) correspond to the minimum and maximum eigenvalues of \(\var^{\policy}_{t}\) and \(\var^{\cpolicy}_{t}\), respectively.  The inequality \((\star)\) follows from \cite[Theorem 4.2.2]{ref:HJ-20}. Let us show that the product \[\prod_{t=1}^n \frac{\det\bigl(\var^{\cpolicy}_{t}\bigr)^{1/2}}{\det\bigl(\var^{\policy}_{t}\bigr)^{1/2}} \] is bounded across \(n\). {\color{black}Indeed, we observe that 
\begin{equation*}
    \begin{aligned}
        \det\bigl(\var^{\cpolicy}_{t}\bigr) &= \det\Bigl(\wvar_1 - \wvar_{12}\bigl(B_H \var_2 B_H^{\top} + \wvar_2 \bigr){\inverse}\wvar_{21}\Bigr) \\& \hspace{35mm}  \det\bigl(B_{H} \var_2 B_{H}^{\top} + \wvar_2 \bigr) \\
        &< \det\Big(B_{M} \var_1 B_{M}^{\top}+ \wvar_1  - \wvar_{12}\bigl(B_H \var_2 B_H^{\top} + \wvar_2 \bigr){\inverse}\wvar_{21} \Big)\\& \hspace{15mm}\det \bigl(B_{H} \var_2 B_{H}^{\top} + \wvar_2 \bigr)\\
        &= \det\bigl(\var^{\policy}_{t}\bigr),
    \end{aligned}
\end{equation*}}
which implies that for every \(n \in \Nz\) we have
\[
    0 < \frac{\det\bigl(\var^{\cpolicy}_{t}\bigr)^{1/2}}{\det\bigl(\var^{\policy}_{t}\bigr)^{1/2}}  < 1,
\]
which means 
\[
0 < \prod_{t=1}^n \frac{\det\bigl(\var^{\cpolicy}_{t}\bigr)^{1/2}}{\det\bigl(\var^{\policy}_{t}\bigr)^{1/2}}  < 1.
\]
Note that the random process \((\rnd{n}{})_{n\in\Nz}\) (see \eqref{eq:RND e}) is a non-negative martingale, and therefore by the martingale convergence theorem \cite[Corollary 3, p. 149, vol. 2]{ref:Shi-19} the limit \(\lim_n \rnd{n}{} \eqqcolon \rnd{\infty}{}\) exists and is finite. 

To see \ref{it:fa2}, suppose that the hypothesis holds, then
\begin{align*}
    &\rnd{n}{(\stt{}{0},\ldots,\stt{}{n})} \\
    & \leqslant \limsup_{n\ra +\infty} \rnd{n}{(\stt{}{0},\ldots,\stt{}{n})}\\
    &\leqslant  \limsup_{n \in \Nz}\prod_{t=1}^n \frac{\det\bigl(\var^{\cpolicy}_{t}\bigr)^{1/2}}{\det\bigl(\var^{\policy}_{t}\bigr)^{1/2}}  \exp{\Biggl[-\frac{1}{2}\Bigl(\rate_n - 1\Bigr) \sum_{t=1}^n  \breve{s_t}\Biggr]}\\
    & <\limsup_{n \in \Nz} \exp{\Biggl[-\frac{1}{2}\Bigl(\rate_n - 1\Bigr) \sum_{t=1}^n  \breve{s_t}\Biggr]}\\
    &\xrightarrow[]{n\to+\infty} 1 \quad  \PP^{\cpolicy}_{\idist}-\text{almost every}\;\sample,
\end{align*}
which implies that \(\rnd{n}{(\stt{}{0},\ldots,\stt{}{n})} < 1\). Since \(\rnd{n}{} \xrightarrow[n\to+\infty]{\PP^{\cpolicy}_{\idist}-\text{a.s.}} \rnd{\infty}{}\), invoking Lebesgue's dominated convergence theorem \cite[Theorem 3, p. 224, vol. 1]{ref:Shi-16} we have \(\EE \expecof[\big]{\rnd{\infty}{}} = \lim_{n \ra +\infty} \EE \expecof[\big]{\rnd{n}{}} = 1\).\footnote{See Lemma \ref{appen:rnd_mean} in Appendix \ref{appen:prob basics} for the proof.} From \cite[Theorem 2, p. 168, vol. 2]{ref:Shi-19} we conclude that the measures \(\PP^{\policy}_{\idist}\) and \(\PP^{\cpolicy}_{\idist}\) are not mutually singular, and from Theorem \ref{p:key} (the converse statement) we assert that a separator in the sense of Definition \ref{def:separator} does not exist. 

For the proof of \ref{it:faaaa2}, we stipulate that the hypotheses hold. Observe that 
\begin{itemize}[leftmargin=*]
    \item \(\sum_{t=1}^n  \breve{s_t} \leqslant q_n< +\infty\) for every \(n \in \Nz\), where \((q_n)_{n \in \Nz}\) is a bounded sequence. Moreover, \(q_{n+1} - q_n \geq \sum_{t=1}^{n+1}  \breve{s_t} - \sum_{t=1}^n  \breve{s_t} > 0\) for all \(n \in \Nz\) implying that \((q_n)_{n \in \Nz}\) is a strictly increasing sequence. From Monotone convergence theorem \cite[Theorem 3.14]{ref:WR-64}, \((q_n)_{n \in \Nz}\) has a limit and it converges to \(q \Let \sup_{n \in \Nz}q_n < +\infty\). This implies that the sequence \(\Bigl(\sum_{t=1}^n  \breve{s_t} \Bigr)_{n \in \Nz}\) is bounded and \(\sum_{t=1}^n  \breve{s_t} \leqslant \sup_{n \in \Nz} q_n \eqqcolon q\) for all \(n \in \Nz\);
    
    \item we have \(0 < \rate_n \leqslant M\) for all \(n \in \Nz\), which implies that \(-\frac{1}{2}(M-1)\sum_{t=1}^n  \breve{s_t} \leqslant -\frac{1}{2}(\rate_n-1)\sum_{t=1}^n  \breve{s_t} < \frac{1}{2}\sum_{t=1}^n  \breve{s_t} \quad \text{for all }n \in \Nz.\)
\end{itemize}
This yields that \(\exp{\bigl[-\frac{1}{2} (\rate_n-1)\sum_{t=1}^n  \breve{s_t}\bigr]} \leqslant \exp{\bigl[\frac{1}{2}\sum_{t=1}^n  \breve{s_t} \bigr]}\). Therefore, 
\begin{align*}
    &\rnd{n}{(\stt{}{0},\ldots,\stt{}{n})}  \\
    & \leqslant \limsup_{n\ra +\infty} \prod_{t=1}^n \frac{\det\bigl(\var^{\cpolicy}_{t}\bigr)^{1/2}}{\det\bigl(\var^{\policy}_{t}\bigr)^{1/2}}  \exp{\Biggl[-\frac{1}{2}\Bigl(\rate_n - 1\Bigr) \sum_{t=1}^n  \breve{s_t}\Biggr]}\\
    &\leqslant \limsup_{n\ra +\infty} \exp{\Biggl[\frac{1}{2}  \sum_{t=1}^n  \breve{s_t}\Biggr]}= \exp{\Biggl[\frac{1}{2}  \limsup_{n\ra +\infty} \sum_{t=1}^n  \breve{s_t}\Biggr]}\\
    &= \exp{\Bigl[\frac{1}{2}  q\Bigr]}
\end{align*}
Since \(\rnd{n}{} \xrightarrow[n\to+\infty]{\PP^{\cpolicy}_{\idist}-\text{a.s.}} \rnd{\infty}{}\), invoking Lebesgue's dominated convergence theorem \cite[Theorem 3, p. 224, vol. 1]{ref:Shi-16} we have \(\EE \expecof[\big]{\rnd{\infty}{}} = \lim_{n \ra +\infty} \EE \expecof[\big]{\rnd{n}{}} = 1\). Arguments similar to \ref{it:fa2} hold at this stage and we assert that the measures \(\PP^{\policy}_{\idist}\) and \(\PP^{\cpolicy}_{\idist}\) are not mutually singular. From Theorem \ref{p:key} we conclude that a separator in the sense of Definition \ref{def:separator} does not exist.

We proceed to \ref{it:ta2}. Assume that the hypotheses hold. Then we write
\begin{align*}
    0 &\leqslant\rnd{\infty}{(\stt{}{0},\ldots,\stt{}{n})} \\&= \lim_{n \ra +\infty} \rnd{n}{(\stt{}{0},\ldots,\stt{}{n})}\\ &=\liminf_{n \ra +\infty} \rnd{n}{(\stt{}{0},\ldots,\stt{}{n})}\\
    & \leqslant  \liminf_{n \in \Nz}\prod_{t=1}^n \frac{\det\bigl(\var^{\cpolicy}_{t}\bigr)^{1/2}}{\det\bigl(\var^{\policy}_{t}\bigr)^{1/2}}  \exp{\Biggl[-\frac{1}{2}\Bigl(\rate_n - 1\Bigr) \sum_{t=1}^n  \breve{s_t}\Biggr]}\\
    & \leqslant \liminf_{n \ra +\infty} \exp{\Biggl[-\frac{1}{2}\Bigl(\rate_n - 1\Bigr) \sum_{t=1}^n  \breve{s_t}\Biggr]}\\
    &\xrightarrow[]{n\to+\infty} 0 \quad  \PP^{\cpolicy}_{\idist}-\text{almost every}\;\sample,
\end{align*}
which implies that \(\rnd{n}{} \xrightarrow[n\to+\infty]{\PP^{\cpolicy}_{\idist}-\text{a.s.}} 0\). Invoking Theorem \ref{p:key} we conclude that there exists a separator in the sense of Definition \ref{def:separator}. This completes the proof of Theorem \ref{prop:main}. \hfill \(\blacksquare\)
\end{pf}}
{\color{black}\begin{rem}
Observe that in the proof of Theorem \ref{prop:main}, the statistics of the private excitation is chosen in such a way that the ratio is bounded; in particular, \[\frac{\det \bigl( \var^{\cpolicy}_{t} \bigr)}{\det \bigl( \var^{\policy}_{t} \bigr)} \leqslant 1 \quad \text{for every }t \in \Nz,  \]  can be ascertained, which plays a crucial role in the rest of the proof. For the specific case of CPSs under FDI attacks considered in Remark \ref{rem:on fdi}, the conditional variances given \(\cphist_t\) and \(\phist_t\) are identical, i.e., \(\var^{\cpolicy}_{t} = \var^{\policy}_{t}\) and \(\det \bigl( \var^{\cpolicy}_{t} \bigr) = \det \bigl( \var^{\policy}_{t} \bigr)\) for each \(t\), i.e.,
\[\frac{\det \bigl( \var^{\cpolicy}_{t} \bigr)}{\det \bigl( \var^{\policy}_{t} \bigr)} = 1 \quad \text{for every }t \in \Nz.\]
However, \(\meantr{\cpolicy}{t}\) and \(\meantr{\policy}{t}\) evolve in a different manner and they capture the essence of the FDI attack over time \(t \in \Nz\).
\end{rem}}

{\color{black}\subsection{Proof of Corollary \ref{prop:hmain}} \label{appen:cor_proof}
First, we consider the case when \(\policy\) is employed under no-attack conditions. The closed loop process under the control action \(\control_t\) is given by 
    \(
        \proc_{t+1} = A \proc_t + B \control_t + W_t\; \text{where }\control_t \sim \policy_t \in \pclass{RM}.
    \)
By definition, \(\policy_t\bigl(\hist_t;\cdot\bigr) = \policy_t(\stt{}{t};\cdot)\) for every \(\hist_t.\)
    For a fixed \(S \in \Borelsigalg(\stsp)\), we have
    \begin{align*}
        \PP^{\policy}_{\idist} \cprobof[\Big]{\proc_{t+1} \in S \given \hist_t} &= \int_\admact \tRR^{\policy}_{\idist} \cprobof[\Big]{\proc_{t+1} \in S, \control_t \in \odif{\cont{}{t}} \given \hist_t}\\
        & = \int_{\admact} \tRR^{\policy}_{\idist} \cprobof[\Big]{\proc_{t+1} \in S \given \hist_t,\cont{}{t}} \,\policy_t(\hist_t;\odif{\cont{}{t}})\\
        & =\int_{\admact} \tRR^{\policy}_{\idist} \cprobof[\Big]{\proc_{t+1} \in S \given \stt{}{t},\cont{}{t}} \,\policy_t(\stt{}{t};\odif{\cont{}{t}})\\
        &=\int_\admact \tRR^{\policy}_{\idist} \cprobof[\Big]{\proc_{t+1} \in S, \control_t \in \odif{\cont{}{t}} \given \stt{}{t}}\\
        &= \PP^{\policy}_{\idist} \cprobof[\Big]{\proc_{t+1} \in S \given \stt{}{t}}.
    \end{align*}
A similar inference can be drawn when \(\cpolicy\) is employed under attack conditions. This completes the proof of the first assertion \ref{it:cor-markov}.
    
The proofs of the assertions \ref{it:cor-gaussian} and \ref{it:cor-th2}, proceed along the same lines as those of Lemma \ref{lem:rvg} and Theorem \ref{prop:main}, respectively, and have been omitted. \hfill \(\blacksquare\)
}

\section{Background on Probability theory}
\label{appen:prob basics}
Consider two probability measures \(\PP, \PP'\) on the same measurable space \((\Omega, \sigalg)\).
\begin{itemize}[label=\(\circ\), leftmargin=*]
	\item Recall that absolutely continuity \cite[p.\ 233, p.\ 438]{ref:Shi-16} of \(\PP\) relative to \(\PP'\) is denoted by \(\PP\abscont\PP'\). If \(\PP\abscont\PP'\), then the \emph{Radon-Nikodym theorem} \cite[p.\ 233]{ref:Shi-16} asserts the existence of a non-negative \emph{density} \(\odv{\PP}{\PP'}:\Omega\lra\lcro{0}{+\infty}\), a Borel measurable and \(\PP'\)-integrable function on \(\Omega\), such that \(\int_A \PP(\odif y) = \int_A \odv{\PP}{\PP'}(y)\:\PP'(\odif{y})\) for all \(A\in\sigalg\).
	\begin{lem}
    \label{appen:rnd_mean}
    \(\EE \expecof[\big]{\odv{\PP}{\PP'}} = 1\) where expectation is taken with respect to the measure \(\PP'\).
    \end{lem}
    \begin{pf}
    Indeed, note that 
    \[\EE \expecof[\Bigg]{\odv{\PP}{\PP'}}=\int_{\Omega} \odv{\PP}{\PP'}(y)\:\PP'(\odif{y}) = \int_{\tss} \PP(\odif y) = 1. 
    \] \hfill \(\blacksquare\)
\end{pf}
	\item \emph{Support} of a probability measure \(\PP\) is defined by the closure of the set \( \aset[\big]{E \in \sigalg \suchthat \PP(E) \neq 0}\).
	
	\item If there exists \(B\in\sigalg\) such that \(\PP(B) = 1 = \PP'(\Omega\setmin B)\), then \(\PP\) and \(\PP'\) are singular \cite[p.\ 276, p.\ 438]{ref:Shi-16}, and we write \(\PP\perp\PP'\). 
\end{itemize}

\begin{defn}[\cite{ref:HLJL-12}]\label{def:sk}
Let \(E_1\) and \(E_2\) be two Borel subsets of \(\Omega\). A stochastic kernel on \(E_1\) given \(E_2\) is a function \(\varphi(\cdot \mid \cdot) \) such that 
\begin{enumerate}
    \item for each \(x \in E_2\), \(\varphi(\cdot \mid y)\) is a probability measure on \(E_1\);
    \item for each \(T \in \Borel(E_1)\), \(\varphi(T\mid \cdot)\) is a Borel measurable function on \(E_2\).
\end{enumerate}
\end{defn}
\begin{defn}
 \label{def:policy}
 An \emph{admissible control policy} \cite[p. 285]{ref:AraBorFerGhoMar-93} is a sequence \(\policy = (\policy_t)_{t \in \Nz}\) of stochastic kernels on the set \(\admact\) given \(\history_t\) that satisfies
 \(
	\policy_t(\history_t; \admact) = 1. 
\)
\end{defn}

\bibliographystyle{plain}
\bibliography{ref}

\begin{thebibliography}{10}

\bibitem{ref:AraBorFerGhoMar-93}
A.~Arapostathis, V.~S. Borkar, E.~Fern\'{a}ndez-Gaucherand, M.~K. Ghosh, and
  S.~I. Marcus.
\newblock Discrete-time controlled {M}arkov processes with average cost
  criterion: a survey.
\newblock {\em SIAM Journal on Control and Optimization}, 31(2):282--344, 1993.

\bibitem{ref:VSB-08}
V.~S. Borkar.
\newblock {\em Stochastic {A}pproximation}.
\newblock Cambridge University Press, Cambridge; Hindustan Book Agency, New
  Delhi, 2008.
\newblock A dynamical systems viewpoint.

\bibitem{ref:DD-GLH-YX-etal-18}
D.~Ding, Q.-L. Han, Y.~Xiang, X.~Ge, and X.-M. Zhang.
\newblock A survey on security control and attack detection for industrial
  cyber-physical systems.
\newblock {\em Neurocomputing}, 275:1674--1683, 2018.

\bibitem{ref:CF-YQ-PC-WXZ-17}
C.~Fang, Y.~Qi, P.~Cheng, and W.~X. Zheng.
\newblock Cost-effective watermark based detector for replay attacks on
  cyber-physical systems.
\newblock In {\em 2017 11th Asian Control Conference (ASCC)}, pages 940--945.
  IEEE, 2017.

\bibitem{ref:AF-WS-18}
A.~Ferdowsi and W.~Saad.
\newblock Deep learning-based dynamic watermarking for secure signal
  authentication in the internet of things.
\newblock In {\em 2018 IEEE International Conference on Communications (ICC)},
  pages 1--6. IEEE, 2018.

\bibitem{ref:HF-67}
H.~Furstenberg.
\newblock Disjointness in ergodic theory, minimal sets, and a problem in
  {D}iophantine approximation.
\newblock {\em Math. Systems Theory}, 1:1--49, 1967.

\bibitem{ref:JG-DU-AC-etal-18}
J.~Giraldo, D.~Urbina, A.~Cardenas, J.~Valente, M.~Faisal, J.~Ruths, N.~O.
  Tippenhauer, H.~Sandberg, and R.~Candell.
\newblock A survey of physics-based attack detection in cyber-physical systems.
\newblock {\em ACM Computing Surveys (CSUR)}, 51(4):1--36, 2018.

\bibitem{ref:RG-CS-EN-SR-22}
R.~Goyal, C.~Somarakis, E.~Noorani, and S.~Rane.
\newblock Co-design of watermarking and robust control for security in
  cyber-physical systems.
\newblock {\em arXiv preprint arXiv:2209.06267}, 2022.

\bibitem{ref:HLJL-12}
O.~H.-Lerma and J.~B. Lasserre.
\newblock {\em Discrete-{T}ime {M}arkov {C}ontrol {P}rocesses: {B}asic
  {O}ptimality {C}riteria}, volume~30.
\newblock Springer Science \& Business Media, 2012.

\bibitem{ref:SH-MX-HHC-YL-14}
S.~Han, M.~Xie, H.-H. Chen, and Y.~Ling.
\newblock Intrusion detection in cyber-physical systems: Techniques and
  challenges.
\newblock {\em IEEE systems journal}, 8(4):1052--1062, 2014.

\bibitem{ref:HH-JY:16}
H.~He and J.~Yan.
\newblock Cyber-physical attacks and defences in the smart grid: a survey.
\newblock {\em IET Cyber-Physical Systems: Theory \& Applications},
  1(1):13--27, 2016.

\bibitem{ref:JRH-LDC-JGA-17}
J.~R. Hernan, L.~D. Cicco, and J.~G. Alfaro.
\newblock On the use of watermark-based schemes to detect cyber-physical
  attacks.
\newblock {\em EURASIP Journal on Information Security}, 2017(1):1--25, 2017.

\bibitem{ref:PH-MP-RV-AA-17}
P.~Hespanhol, M.~Porter, R.~Vasudevan, and A.~Aswani.
\newblock Dynamic watermarking for general {LTI} systems.
\newblock In {\em 2017 IEEE 56th Annual Conference on Decision and Control
  (CDC)}, pages 1834--1839. IEEE, 2017.

\bibitem{ref:PH-MP-RV-AA-18}
P.~Hespanhol, M.~Porter, R.~Vasudevan, and A.~Aswani.
\newblock Statistical watermarking for networked control systems.
\newblock In {\em 2018 Annual American Control Conference (ACC)}, pages
  5467--5472. IEEE, 2018.

\bibitem{ref:PH-MP-RV-AA-20}
P.~Hespanhol, M.~Porter, R.~Vasudevan, and A.~Aswani.
\newblock Sensor switching control under attacks detectable by finite sample
  dynamic watermarking tests.
\newblock {\em IEEE Transactions on Automatic Control}, 66(10):4560--4574,
  2020.

\bibitem{ref:HJ-20}
R.~A. Horn and C.~R. Johnson.
\newblock {\em Matrix {A}nalysis}.
\newblock Cambridge University Press, 2012.

\bibitem{ref:MH-TT-VG-16}
M.~Hosseini, T.~Tanaka, and V.~Gupta.
\newblock Designing optimal watermark signal for a stealthy attacker.
\newblock In {\em 2016 European Control Conference (ECC)}, pages 2258--2262.
  IEEE, 2016.

\bibitem{ref:TH-BS-PRK-LX-18}
T.~Huang, B.~Satchidanandan, P.~R. Kumar, and L.~Xie.
\newblock An online detection framework for cyber attacks on automatic
  generation control.
\newblock {\em IEEE Transactions on Power Systems}, 33(6):6816--6827, 2018.

\bibitem{ref:AH-JL-FL-BL-17}
A.~Humayed, J.~Lin, F.~Li, and B.~Luo.
\newblock Cyber-physical systems security—a survey.
\newblock {\em IEEE Internet of Things Journal}, 4(6):1802--1831, 2017.

\bibitem{ref:MJK-AK-MF-TJ-19}
M.~J. Khojasteh, A.~Khina, M.~Franceschetti, and T.~Javidi.
\newblock Authentication of cyber-physical systems under learning-based
  attacks.
\newblock {\em IFAC-PapersOnLine}, 52(20):369--374, 2019.

\bibitem{ref:MJK-AK-MF-TJ-20}
M.~J. Khojasteh, A.~Khina, M.~Franceschetti, and T.~Javidi.
\newblock Learning-based attacks in cyber-physical systems.
\newblock {\em IEEE Transactions on Control of Network Systems}, 8(1):437--449,
  2020.

\bibitem{ref:WHK-BS-PRK-19}
W.-H. Ko, B.~Satchidanandan, and P.~R. Kumar.
\newblock Dynamic watermarking-based defense of transportation cyber-physical
  systems.
\newblock {\em ACM Transactions on Cyber-Physical Systems}, 4(1):1--21, 2019.

\bibitem{ref:DK-13}
D.~Kushner.
\newblock The real story of stuxnet.
\newblock {\em IEEE Spectrum}, 50(3):48--53, 2013.

\bibitem{ref:LevPelPer-15}
N.~Lev, R.~Peled, and Y.~Peres.
\newblock {Separating signal from noise}.
\newblock {\em Proceedings of the London Mathematical Society},
  110(4):883--931, 2015.

\bibitem{ref:KL-YM-KHJ-21}
H.~Liu, Y.~Mo, and K.~H. Johansson.
\newblock Active detection against replay attack: A survey on watermark design
  for cyber-physical systems.
\newblock In {\em Safety, Security and Privacy for Cyber-Physical Systems},
  pages 145--171. Springer, 2021.

\bibitem{ref:HL-JY-YM-KHJ-18}
H.~Liu, J.~Yan, Y.~Mo, and K.~H. Johansson.
\newblock An on-line design of physical watermarks.
\newblock In {\em 2018 IEEE Conference on Decision and Control (CDC)}, pages
  440--445. IEEE, 2018.

\bibitem{ref:RM-IRC-14}
R.~Mitchell and I.-R. Chen.
\newblock A survey of intrusion detection techniques for cyber-physical
  systems.
\newblock {\em ACM Computing Surveys (CSUR)}, 46(4):1--29, 2014.

\bibitem{ref:YM-RC-BS-13}
Y.~Mo, R.~Chabukswar, and B.~Sinopoli.
\newblock Detecting integrity attacks on scada systems.
\newblock {\em IEEE Transactions on Control Systems Technology},
  22(4):1396--1407, 2013.

\bibitem{ref:YM-EG-AC-BS-10}
Y.~Mo, E.~Garone, A.~Casavola, and B.~Sinopoli.
\newblock False data injection attacks against state estimation in wireless
  sensor networks.
\newblock In {\em 49th IEEE Conference on Decision and Control (CDC)}, pages
  5967--5972. IEEE, 2010.

\bibitem{ref:YM-BS-09}
Y.~Mo and B.~Sinopoli.
\newblock Secure control against replay attacks.
\newblock In {\em 2009 47th annual Allerton conference on communication,
  control, and computing (Allerton)}, pages 911--918. IEEE, 2009.

\bibitem{ref:YM-SW-BS-15}
Y.~Mo, S.~Weerakkody, and B.~Sinopoli.
\newblock Physical authentication of control systems: Designing watermarked
  control inputs to detect counterfeit sensor outputs.
\newblock {\em IEEE Control Systems Magazine}, 35(1):93--109, 2015.

\bibitem{ref:MO-SS-PH-MP-RV-AA-20}
M.~Olfat, S.~Sloan, P.~Hespanhol, M.~Porter, R.~Vasudevan, and A.~Aswani.
\newblock Covariance-robust dynamic watermarking.
\newblock In {\em 2020 59th IEEE Conference on Decision and Control (CDC)},
  pages 3793--3799. IEEE, 2020.

\bibitem{ref:MP-SD-AJ-PH-AA-MJR-RV-20}
M.~Porter, S.~Dey, A.~Joshi, P.~Hespanhol, Anil A.~Aswani, Matthew
  Johnson-Roberson, and Ram Vasudevan.
\newblock Detecting deception attacks on autonomous vehicles via linear
  time-varying dynamic watermarking.
\newblock In {\em 2020 IEEE Conference on Control Technology and Applications
  (CCTA)}, pages 1--8. IEEE, 2020.

\bibitem{ref:MP-PH-AA-MJR-RV-20}
M.~Porter, P.~Hespanhol, A.~Aswani, M.~J. Roberson, and R.~Vasudevan.
\newblock Detecting generalized replay attacks via time-varying dynamic
  watermarking.
\newblock {\em IEEE Transactions on Automatic Control}, 66(8):3502--3517, 2020.

\bibitem{ref:MP-AJ-SD-QW-PH-AA-MJR-RV}
M.~Porter, A.~Joshi, S.~Dey, Q.~Wu, P.~Hespanhol, A.~Aswani, M.J. Roberson, and
  R.~Vasudevan.
\newblock Resilient control of platooning networked robotic systems via dynamic
  watermarking.
\newblock {\em arXiv preprint arXiv:2106.07541}, 2021.

\bibitem{ref:MLP-14}
M.~L. Puterman.
\newblock {\em Markov {D}ecision {P}rocesses: {D}iscrete {S}tochastic {D}ynamic
  {P}rogramming}.
\newblock John Wiley \& Sons, 2014.

\bibitem{ref:MAR-AGA-13}
M.~A. Rahimian and A.~G. Aghdam.
\newblock Structural controllability of multi-agent networks: Robustness
  against simultaneous failures.
\newblock {\em Automatica}, 49(11):3149--3157, 2013.

\bibitem{ref:HS-SA-KHJ-15}
H.~Sandberg, S.~Amin, and K.~H. Johansson.
\newblock Cyberphysical security in networked control systems: An introduction
  to the issue.
\newblock {\em IEEE Control Systems Magazine}, 35(1):20--23, 2015.

\bibitem{ref:HS-VG-KHJ:22}
H.~Sandberg, V.~Gupta, and K.~H. Johansson.
\newblock Secure networked control systems.
\newblock {\em Annual Review of Control, Robotics, and Autonomous Systems},
  5:445--464, 2022.

\bibitem{ref:BS-PRK-16}
B.~Satchidanandan and P.~R. Kumar.
\newblock Dynamic watermarking: Active defense of networked cyber--physical
  systems.
\newblock {\em Proceedings of the IEEE}, 105(2):219--240, 2016.

\bibitem{ref:BS-PRK-CDC-16}
B.~Satchidanandan and P.~R. Kumar.
\newblock Secure control of networked cyber-physical systems.
\newblock In {\em 2016 IEEE 55th Conference on Decision and Control (CDC)},
  pages 283--289. IEEE, 2016.

\bibitem{ref:BS-PRK-17}
B.~Satchidanandan and P.~R. Kumar.
\newblock On minimal tests of sensor veracity for dynamic watermarking-based
  defense of cyber-physical systems.
\newblock In {\em 2017 9th International Conference on Communication Systems
  and Networks (COMSNETS)}, pages 23--30. IEEE, 2017.

\bibitem{ref:BS-PRK-19}
B.~Satchidanandan and P.~R. Kumar.
\newblock On the design of security-guaranteeing dynamic watermarks.
\newblock {\em IEEE Control Systems Letters}, 4(2):307--312, 2019.

\bibitem{ref:NS-SM-17}
N.~Sayfayn and S.~Madnick.
\newblock Cybersafety analysis of the maroochy shire sewage spill (preliminary
  draft).
\newblock {\em MIT Press}, pages 1--29, 2017.

\bibitem{ref:Shi-16}
A.~N. Shiryaev.
\newblock {\em Probability. 1}, volume~95 of {\em Graduate Texts in
  Mathematics}.
\newblock Springer, New York, 3rd edition, 2016.
\newblock Translated from the fourth (2007) Russian edition by R. P. Boas and
  D. M. Chibisov.

\bibitem{ref:Shi-19}
A.~N. Shiryaev.
\newblock {\em Probability. 2}, volume~95 of {\em Graduate Texts in
  Mathematics}.
\newblock Springer, New York, 3rd edition, 2019.
\newblock Translated from the 2007 fourth Russian edition by R. P. Boas and D.
  M. Chibisov.

\bibitem{ref:Str-11}
D.~W. Stroock.
\newblock {\em Probability {T}heory}.
\newblock Cambridge University Press, Cambridge, 2nd edition, 2011.
\newblock An analytic view.

\bibitem{ref:ST-DD-WZS-JY-SKD-16}
S.~Tan, D.~De, W.-Z. Song, J.~Yang, and S.~K. Das.
\newblock Survey of security advances in smart grid: a data driven approach.
\newblock {\em IEEE Communications Surveys \& Tutorials}, 19(1):397--422, 2016.

\bibitem{ref:JT-JS-AG-21}
J.~Tang, J.~Song, and A.~Gupta.
\newblock A dynamic watermarking algorithm for finite markov decision problems.
\newblock {\em arXiv preprint arXiv:2111.04952}, 2021.

\bibitem{ref:PV-13}
P.~Venkitasubramaniam.
\newblock Privacy in stochastic control: a markov decision process perspective.
\newblock In {\em 2013 51st Annual Allerton Conference on Communication,
  Control, and Computing (Allerton)}, pages 381--388. IEEE, 2013.

\bibitem{ref:WR-64}
Walter W.~Rudin.
\newblock {\em Principles of {M}athematical {A}nalysis}.
\newblock McGraw-Hill Book Co., New York, second edition, 1964.

\bibitem{ref:SW-YM-BS-14}
S.~Weerakkody, Y.~Mo, and B.~Sinopoli.
\newblock Detecting integrity attacks on control systems using robust physical
  watermarking.
\newblock In {\em 53rd IEEE Conference on Decision and Control}, pages
  3757--3764. IEEE, 2014.

\bibitem{ref:KZ-16}
K.~Zetter.
\newblock Inside the cunning, unprecedented hack of {U}kraine’s power grid.
\newblock {\em Wired}, 9:1--5, 2016.

\bibitem{ref:LZ-KGV-JH-21}
L.~Zhai, K.~G. Vamvoudakis, and J.~Hugues.
\newblock Switching watermarking-based detection scheme against replay attacks.
\newblock In {\em 2021 60th IEEE Conference on Decision and Control (CDC)},
  pages 4200--4205. IEEE, 2021.

\bibitem{ref:JZ-LP-QLH-CC-SW-YX-21}
J.~Zhang, L.~Pan, Q.-L. Han, C.~Chen, S.~Wen, and Y.~Xiang.
\newblock Deep learning based attack detection for cyber-physical system
  cybersecurity: A survey.
\newblock {\em IEEE/CAA Journal of Automatica Sinica}, 9(3):377--391, 2021.

\end{thebibliography}

\end{document}